\documentclass[10pt,twosided]{article}
\usepackage[tbtags]{amsmath}
\usepackage{epsfig,amstext,amssymb,amsthm,latexsym}
\allowdisplaybreaks[4]
\usepackage{color}

\usepackage{amssymb,color}

\definecolor{c20}{rgb}{0.,0.7,0.}
\definecolor{c30}{rgb}{0.,0.,1.}
\definecolor{c40}{rgb}{1,0.1,0.7}
\definecolor{c50}{rgb}{1,0,0}
\definecolor{c60}{rgb}{1,0.9,0.1}

\def\Hr#1{\textcolor{c20}{#1}}
\def\Hr#1{#1}

\def\aE#1{\textcolor{c20}{#1}}
\def\aE#1{#1}

\def\cE#1{\textcolor{c20}{#1}}
\def\cE#1{#1}

\def\cH#1{\textcolor{c20}{#1}}
\def\cH#1{#1}

\def\cL#1{#1}

\def\cJ#1{\textcolor{c30}{#1}}
\def\cJ#1{#1}

\usepackage{amssymb,color}

\newcommand{\kb}[1]{\boldsymbol{#1}}
\newcommand{\vk}[1]{\kb{#1}}
\def\kal#1{{\cal{ #1}}}

\newcommand{\ve}{\varepsilon}

\newcommand{\abs}[1]{\lvert #1 \rvert}

\newcommand{\E}[1]{\mathbb{E}\left(#1\right)}

\newcommand{\eg}[1]{\mathbb{E}\biggl\{#1 \biggr\}}

\newcommand{\pk}[1]{\mathbb{P} \left( #1 \right) }
\newcommand{\expo}[1]{\exp \left( #1 \right) }

\newcommand{\R}{\!I\!\!R}

\newcommand{\kk}{\!I\!\!k}
\newcommand{\ii}{\!I\!\!I_{\kk}}
\newcommand{\iii}{{\!I\!\!I}'_{\kk'}}
\newcommand{\N}{\!I\!\!N}

\newcommand{\inn}{\in \N}

\newcommand{\ldot}{,\ldots,}

\newcommand{\LL}[1]{\mathcal{L}_{#1} }
\newcommand{\UU}[1]{\mathcal{U}_{#1} }

\newcommand{\BQN}{\begin{eqnarray}}
\newcommand{\EQN}{\end{eqnarray}}
\newcommand{\BQNY}{\begin{eqnarray*}}
\newcommand{\EQNY}{\end{eqnarray*}}

\newcommand{\BS}{\begin{sat}}
\newcommand{\ES}{\end{sat}}
\newcommand{\BT}{\begin{theo}}
\newcommand{\ET}{\end{theo}}
\newcommand{\BK}{\begin{korr}}
\newcommand{\EK}{\end{korr}}

\newcommand{\BD}{\begin{de}}
\newcommand{\ED}{\end{de}}

\newcommand{\BIT}{\begin{itemize}}
\newcommand{\EIT}{\end{itemize}}
\newcommand{\BDI}{\begin{description}}
\newcommand{\EDI}{\end{description}}

\newcommand{\BRM}{\begin{remark}}
\newcommand{\ERM}{\end{remark}}

\newcommand{\BEL}{\begin{lem}}
\newcommand{\EEL}{\end{lem}}

\newcommand{\BL}{\begin{lem}}
\newcommand{\EL}{\end{lem}}

\newtheorem{theo}{Theorem}[section]
\newtheorem{sat}[theo]{Proposition}
\newtheorem{de}[theo]{Definition}
\newtheorem{lem}[theo]{Lemma}

\newtheorem{korr}[theo]{Corollary}
\newtheorem{remark}[theo]{Remark}

\newcommand{\nelem}[1]{{Lemma \ref{#1}}}

\newcommand{\netheo}[1]{{Theorem \ref{#1}}}

\newcommand{\prooftheo}[1]{ \textsc{Proof of Theorem} \ref{#1} }

\newcommand{\prooflem}[1]{\textsc{Proof of Lemma} \ref{#1}}

\newcommand{\COM}[1]{}

\newcommand{\QED}{\hfill $\Box$}

\def\BL{\cE{\beta}}

\def\rw{\rightarrow}
\newcommand{\mq}[1]{\quad\!\!\!\! \mathrm{#1}\quad\!\!\!\!}

\newcommand{\tf}[1]{\vk{#1}}

\def\rw{\rightarrow}
\def\vn{\varepsilon}
\def\mH{\mathcal{H}}
\def\IF{\infty}
\def\gok{1,\cdots,k_1}
\def\gokkk{1,\cdots,k_2}
\def\gokk{k_1+1,\cdots,k}

\def\tLu{ \cE{\eta(u,\cH{k_1}, \vk{\alpha}, \vk{\beta})}}
\def\BB{\mathcal{B}}

\def\AAk{\mathcal{A}_{\kk}}
\def\ovB{\overline{B}}

\def\ZZ{\mathbb{Z}}

\date{}

\def\alpt{\alpha(\tf{t})}

\def\Xjpik{X_{\tf{j,p,},u}^{\ii,\kk}}
\def\Ajpik{A_{\vk{j,p}}^{\ii,\kk}}
\def\Ajpikp{A_{\vk{j',p'}}^{\iii,\kk'}}

\def\ee#1{\textcolor{c30}{#1}}
\def\ee#1{#1}
\def\bb#1{\textcolor{c30}{#1}}
\def\bb#1{#1}
\def\bbb#1{\textcolor{c30}{#1}}
\def\bbb#1{#1}
\def\cL#1{\textcolor{c30}{#1}}
\def\cL#1{#1}

\begin{document}
\COM{Extremes of alpha(t)-locally Stationary Gaussian Random Fields
The main result of this contribution is the derivation of the exact asymptotic behaviour of the supremum of \alpha(t)-locally stationary Gaussian random field over a finite hypercube. We present two applications of our results; the first one deals with the extremes of \cL{aggregate} multifractional Brownian motions, whereas the second one establishes the exact asymptotics of the supremum of \chi-processes generated by multifractional Brownian motions.

\alpha(t)-locally stationary random fields;  Fractional Brownian motion; Multifractional Brownian motion;
Gaussian processes; \chi-processes; Gaussian random fields; metric entropy; weak convergence;  Pickands constant; Piterbarg-Prisyazhnyuk constant
}

\centerline{\bf \Large \ee{Extremes of}  $\alpt$-locally Stationary Gaussian \cE{Random} Fields}

\vskip 1 cm
\centerline{Enkelejd Hashorva and   Lanpeng Ji\footnote{%Department of Actuarial Science, %Faculty of Business and Economics,
University of Lausanne, UNIL-Dorigny 1015 Lausanne, Switzerland\\
%Both authors kindly acknowledge partial financial support by the Swiss National Science Foundation Grant 200021-1401633/1\\
}%; Piterbarg-Prisyazhnyuk constant
}
%\centerline{University of Lausanne, Switzerland}

% \maketitle
\vskip 0.5 cm
        \centerline{\today{}}

%\vskip 1 cm
{\bf Abstract:} The main result of this contribution is the derivation of the exact asymptotic behaviour of the supremum of $\alpha(\vk{t})$-locally stationary \cJ{Gaussian} random field over a finite hypercube. \aE{We present two applications of our results; the} first one deals with the extremes of \cL{aggregate} multifractional \cJ{Brownian motions}, whereas the second \cJ{one} establishes the exact asymptotics of the supremum of $\chi$-processes generated by multifractional Brownian motions.

\section{Introduction and Main Result}
The classical Central Limit Theorem and its ramifications show that \aE{the} Gaussian model is a natural and correct paradigm for building an approximate solution to many otherwise
unsolvable problems encountered in various research fields. While the theory of Gaussian processes and Gaussian random fields (GRF's) is well-developed and mature, the range of \bb{their} applications is constantly growing.   \bb{Recently}, applications in
brain mapping, cosmology, quantum chaos and some other fields have been added to its palmares, see e.g., Adler (2000), Adler and Taylor (2007), Anderes and Chatterjee (2009),  \cL{Aza\"{i}s and Wschebor (2009) and  Adler et al. (2012a,b).} In applications related to extremes of
\Hr{Gaussian processes}  the fractional Brownian motion (fBm) appears inevitably in the definition of the Pickands constant, see e.g., \aE{Pickands (1969)}, Berman (1992) and Piterbarg (1996).
Numerous research articles have shown the importance of fBm in \aE{both} theoretical models and applications. For certain applications, the stationarity of increments, which together with the
self-similarity property characterises  fBm  in the class of Gaussian processes can be a severe  restriction. A natural way to avoid the stationarity of increments property is to introduce the multifractional Brownian motion (mfBm), see e.g., Stoev and Taqqu (2006) and
Ayachea  et al. (2011). \aE{In order to make the problem tractable, we discuss in this paper a simple class of mfBm.} By definition, a mean-zero
Gaussian process  $\{\cJ{B_{\alpha(t)}}(t),t\ge 0\}$ is \cJ{called a} mfBm
%$\{B_{H(t)}(t),t\ge0\}$ be a multifractional Brownian motion (mfBm)
with parameter $\alpha(t),t\ge 0,$ if
\BQN
\E{B_{\alpha(t)}(t) B_{\alpha(s)}(s)}=\frac{1}{2}D(\alpha(s,t))\left(s^{\alpha(s,t)}+t^{\alpha(s,t)}-|t-s|^{\alpha(s,t)}\right), \ \
\alpha(s,t):=\cJ{\alpha(s)/2+\alpha(t)/2}, s,t\ge 0,
\EQN
where $D(x)=\frac{2\pi}{\Gamma(x+1)\sin\left(\pi x/2\right)}$ and
$\alpha(\cdot)$ is a H\"{o}lder function of exponent \bbb{$\gamma>0$} such that $0<\alpha(t)<\cJ{2}\min(1,\gamma),t\ge 0$, see e.g.,  Ayache et al.\ (2000).  For $\alpha(t)=\alpha\in (0,2), t\ge 0,$ the $B_\alpha$ reduces to a fBm (\cJ{not necessarily standard}).\\
Inspired by the structure of the mfBm, the recent  paper D\c{e}bicki and Kisowski (2008) introduces the class of $\alpha(t)$-locally stationary
Gaussian processes.  Therein the exact asymptotic  of the tail behaviour of the supremum of $\alpha(t)$-locally stationary  Gaussian process
is derived which can be applied, for instance, to analyse the extremes of standardized mfBm.\\
\aE{It is worth noting that this new class includes locally stationary ones, see Berman (1974), H\"usler (1990)
 \Hr{and Piterbarg (1996)} for results concerning the \Hr{asymptotic} behaviour of their extremes.}  If $\{X_i(t), t\in [0,T]\},i\le k,$ are independent real-valued Gaussian processes  a natural \cJ{GRF}
associated with these processes is the \cL{aggregate} random field
$$X(\vk{t})= \sum_{i=1}^k X_i(t_i), \quad \vk{t}=(t_1 \ldot t_k)\in [0,T]^k.$$
 Extremes of GRF's can not be analysed by aggregating \cJ{the }corresponding results for processes. Moreover, the analysis of the extremes of  GRF's leads to technical difficulties, see e.g., the excellent monographs Piterbarg (1996) and Adler and Taylor (2007). Recently,
Abramowicz and Seleznjev (2011) deal with multivariate
 piecewise linear interpolation of locally stationary random fields, whereas Hashorva et al.\ (2012) investigates the piece-wise approximation of
 $\alpha(t)$-locally stationary processes. With motivation from  the aforementioned papers and
 D\c{e}bicki and Kisowski (2008), we consider, in this paper, extremes of
$\alpha(\tf{t})$-locally stationary  GRF   $\{X(\tf{t}), \tf{t}\in [0,T]^k\}$ \cL{(to be defined below)}.
 Specifically, we are interested \cJ{in} the exact asymptotic behaviour of
\BQN
\pk{\underset{\tf{t}\in [0,T]^k}{\sup}X(\tf{t})>u}, \quad \quad u\rightarrow\infty,
\label{localsta1}
\EQN
 with $T>0$ a given constant and $k\inn$ a positive integer. %and $X(\tf(t))$ an $\alpha(\tf{t})$

\COM{
We first introduce some notation before giving the definition of the $\alpha(\tf{t})$-locally stationary  GRF's  . Let $\tf{l}=(l_1,\cdots,l_k)$ be a vector of positive integers such that $\sum_{j=1}^k l_j=d,$ and let $L_i=\sum_{j=1}^i l_j, i=1,\cdots,k, L_0=0$ be the sequence of its cumulative sums. Thus, the vector $\tf{l}$ defines the $l$-decomposition of $[0,T]^k$ into $\prod_{j=1}^k [0,T]^{l_j}$. For $\tf{t}\in [0,T]^k$ we denote the coordinates vector corresponding to the $j$-th component of the decomposition by
$$
\tf{t}^j:=(t_{L_{j-1}+1},\cdots,t_{L_j}).
$$
In the following, $||\tf{t}||$ represents the Euclidean norm of vector $\tf{t}$. The pair $(L,\alpha)$ is called a {\it structure}. For any $\tf{t}=(t_1,\cdots,t_d)$ and $\alpha=(\alpha_1,\cdots,\alpha_k)$, the {\it structure modulus} is defined by
\BQNY
|\tf{t}|_{(L,\alpha)}=\sum_{i=1}^k\left(\sum_{j=L_{i-1}+1}^{L_i}(t_j)^2\right)^{\alpha_i/2}
\EQNY

{\bf Note: in the following, the letter $'b'$ in the notation of the form $"x^b"$ is a script rather than a power exponent.}
}

\bbb{Let $\mathcal{C}(\tf{D})$ denote the set \cE{of} %including
all continuous functions on $\tf{D}\subset\R^{\cH{k}}$. Next, we give a formal definition of the GRF's of interest.
}

{\bf Definition.} A real-valued separable GRF  $\{X(\tf{t}), \tf{t}\in [0,T]^k\}$ is said to be $\alpha(\tf{t})$-locally stationary if

$D1.$ $\E{X(\tf{t})}=0$ and $Var(X(\tf{t}))=1$ for all $\tf{t}\in [0,T]^k$;

$D2.$ $\alpha_i(t_i)\in \mathcal{C}([0,T])$ and $\alpha_i(t_i)\in(0,2]$ for all $t_i\in [0,T]$, $i=1,\cdots,k$;

$D3.$ $C_i({\cJ{\tf{t}}})\in\mathcal{C}([0,T]^k)$ and $0<\inf\{C_i({\cJ{\tf{t}}}): \tf{t}\in[0,T]^k\}\le \sup\{C_i({\cJ{\tf{t}}}): \tf{t}\in[0,T]^k\}:=C_U^i<\infty$, $i=1,\cdots,k$;

$D4.$ \cE{uniformly with respect to $\tf{t}\in [0,T]^k$}
\BQN
1-Cov(X(\tf{t}),X(\tf{t}+\tf{s}))=\sum_{i=1}^k C_i({\cJ{\tf{t}}})|s_i|^{\alpha_i(t_i)}+o\left(\sum_{i=1}^k C_i({\cJ{\tf{t}}})|s_i|^{\alpha_i(t_i)}\right)\label{eqcov1}
\EQN
as $\tf{s}\rightarrow \tf{0}$ with  $\tf{0}:=(0,\cdots,0)\in\R^k$.

A canonical example of $\alpha(\vk{t})$-locally stationary GRF's is \cL{the aggregate} \bbb{mfBm} defined by aggregating independent standardized mfBm's, see Section 2.

In this paper we consider the case that there exists some integer $k_1\le k$ such that:

$A1.$ each of $\alpha_i(t_i), i=1,\cdots,k_1,$ attains its global minimum on $[0,T]$ at a unique point $t_i^0$, and \cL{further for any $i=k_1+1,\cdots,k,$ there is some $\cJ{[a_i,b_i]\subset(0,T)}$ such that $\alpha_i(t_i)\equiv\alpha_i$ in $[a_i,b_i]$ which is the global minimum of $\alpha_i(t_i)$ on $[0,T]$;}

%In order to obtain the exact tail asymptotics of the supremum of $\alpha(\vk{t})$-ls GRF's we need to impose some asymptotic conditions on $\alpha_i$'s, namely\\
%A1. each of $\alpha_i(t_i), i=1,\cdots,k_1,$ attains its global minimum at the unique point $t_i^0\in [0,T]$;

$A2.$ there exist $\cJ{M_i}, \beta_i>0,$ and $\delta_i>1,$ $i=1,\cdots,k_1,$ %\cJ{positive increasing functions $g_i(x), i=1,\cdots,k_1,$ satisfying
%\BQNY
%\underset{x\rw\IF}\lim g_i(R_i\left(\frac{\ln\ln x}{\ln x}\right)^{1/\beta_i})\ln\ln x =0, \ \ \forall R_i>0,
%\EQNY}
 such that
\BQN\label{DK8}
\alpha_i(t_i+t_i^0)=\alpha_i(t_i^0)+M_i|t_i|^{\beta_i}+o(|t_i|^{\beta_i}\abs{\ln\abs{t_i}}^{-\delta_i}),  \quad  \mathrm{as}\ t\aE{\rightarrow 0},
\EQN
and  there exist $\cJ{M_i}, \beta_i, \tilde{M}_i, \tilde{\beta}_i>0,$ and $\tilde{\delta}, \delta_i>1,$ $i=k_1+1,\cdots,k,$
 such that
\BQN
\alpha_i(b_i+t_i)&=&\alpha_i(b_i)+M_it_i^{\beta_i}+o(t_i^{\beta_i}\abs{\ln{t_i}}^{-\delta_i}),  \quad  \mathrm{as}\ t_i\aE{\downarrow 0},\\
\alpha_i(a_i-t_i)&=&\alpha_i(a_i)+\tilde{M}_it_i^{\tilde{\beta}_i}+o(t_i^{\tilde{\beta}_i}\abs{\ln{t_i}}^{-\tilde{\delta}_i}),  \quad  \mathrm{as}\ t_i\aE{\downarrow 0}.\label{DK9}
\EQN

\bb{The assumption $A1$ is initially suggested in} D\c{e}bicki and Kisowski (2008), whereas assumption $A2$ is a weaker version of a similar condition given therein which assumes (\ref{DK8}-\ref{DK9}) with $\abs{t_i}^{\delta_i} (t_i^{\tilde{\delta}_i})$ instead of $\abs{\ln\abs{t_i}}^{-\delta_i} (\abs{\ln{t_i}}^{-\tilde{\delta}_i})$.

For notational simplicity, \bb{set}
$$\alpha_i:=\alpha_i(t^0_i), \quad i=1,\cdots,k_1,$$%,
\cJ{ and
\def\x{\vk{x}}
$$ \int_{\vk{x}\in\{\tf{x}_0\}\times\tf{D}_1}C(\vk{x})d\vk{x}:=\int_{{\x}\in\tf{D}_1}C(\tf{x}_0,\x)d{\x}$$
for all integrable function $C(\cdot)$.
Further, %define, for $u>0$,
%$$\Psi(u):=\frac{1}{\sqrt{2\pi}}\int_u^\IF e^{-\frac{x^2}{2}}dx.$$ %, \mq{and}\tLu:=\frac{\prod_{i=1}^{k}u^{2/\alpha_i}}{\prod_{i=1}^{k_1}(\ln u)^{1/\beta_i}},
%with $\prod_{i=m+1}^m(\cdot):=1, m\inn$,
denote by $\Psi(\cdot)$ the survival function a standard normally distributed random variable, and by $\Gamma(\cdot)$ the Euler's Gamma function.}

The crucial step of the proof of our main result \netheo{mainthm}
is an application of the double-sum method \cL{that} was \cL{developed by Pickands (1969)}. As expected, \cJ{the Pickands constant  defined by}
$$\mathcal{H}_{\cJ{\alpha}}=\lim_{\bbb{\mathcal{T}}\rightarrow\infty} \bbb{\mathcal{T}}^{-1} \eg{ \exp\biggl(\sup_{t\in[0,\bbb{\mathcal{T}}]}\Bigl(\sqrt{2}B_{\bbb{\alpha}}(t)-t^{\alpha}\Bigr)\biggr)}\in (0,\IF),\quad
 \cJ{\alpha\in (0,2]},$$
\cJ{appears in the asymptotic expansion}, where $\{B_{\bbb{\alpha}}(t),t\ge 0\}$ is a fBm with Hurst \Hr{index} $\alpha/2.$
See  Pickands (1969), Piterbarg (1996) or D\c{e}bicki (2002) for the basic properties of Pickands constant and generalisations.

\BT\label{mainthm}
Let $\{X(\tf{t}), \tf{t}\in[0,T]^k\}$ be an $\alpt$-locally stationary GRF that satisfies
\BQN
Cov(X(\tf{t}),X(\tf{s}))<1\label{eqcov11}, \quad \forall \tf{t,s}\in[0,T]^k, \quad \tf{t}\neq\tf{s}.
\EQN
If both conditions $A1$ and $A2$ are satisfied, then we have (set $q_{k_1}:=\#\{i\in\N:1\le i\le k_1, t^0_i\in(0,T)\}$)
\BQN \label{EM}
\pk{\underset{\tf{t}\in [0,T]^k}{\sup}X(\tf{t})>u}&=&\kal{K}_{\cH{\vk{O}}}
 \cL{ u^{\alpha}(\ln u)^\beta} \Psi(u)(1+o(1)), \quad u\rw\IF,
\EQN
\cL{where $\alpha=2\sum_{i=1}^k1/\alpha_i$, $\beta=-\sum_{i=1}^{k_1}1/\beta_i$ and }
\BQN
\kal{K}_{\cH{\vk{O}}}
&=& 2^{q_{k_1}}\biggl(\prod_{i=1}^{k_1}%C_i^{1/\alpha_i}
\Bigl( \frac{\alpha_i^2}{2M_i}\Bigr)^{1/\beta_i}  \cH{\Gamma(1/\beta_i+1)} \biggr)
%\frac{\Gamma(1/\beta_i)}{\beta_i}
\Bigr(\prod_{i=1}^k\mH_{\alpha_i} \Bigl)
\int_{\cH{\vk{x}\in \vk{O}}}\cJ{\prod_{i=1}^{k}(C_i(\tf{x}))^{1/\alpha_i}d\tf{x}}
%\left(\prod_{i=k_1+1}^{k}\int_{a_i}^{b_i}(C_i(x))^{1/\alpha_i}dx\right)
\cH{\in (0,\IF)},
%\int_{\R^{k_1}_+}e^{-\sum_{i=1}^{k_1}\frac{2}{\alpha_i^2}x_i^{\beta_i}}d\tf{x}
\EQN
with $\cJ{\vk{O}=\prod_{i=1}^{k_1}\{t^0_i\}\times\prod_{i=k_1+1}^k [a_i,b_i]}$.
%and $\mH_\alpha$ is the Pickands' constant.% and $C_i=C_i(t^0_i), i=\gok$.
%, and
%$\Psi(u)=\frac{1}{\sqrt{2\pi}}\int_u^\IF e^{-\frac{x^2}{2}}dx$.
\ET
%\cH{The next theorem, which is similar in spirit with \netheo{mainthm} is of particular interest for the derivation of the exact asymptotics of the tail of supremum of $\chi$-processes.}

{\bf Remarks}: a) Under the conditions of \netheo{mainthm}, if, for \cL{the chosen} $k_1<k$,
 $\alpha_i(t_i)\equiv\alpha_i, i=k_1+1,\cdots,k,$ on \cL{some compact set} $\vk{O}_2\cE{\subset}\ \R_+^{\cL{k-k_1}}$, with positive Lebesgue measure, then
 \eqref{EM} holds \cJ{for $\{X(\tf{t}), \tf{t}\in[0,T]^{k_1}\times \vk{O}_2\}$} with $\vk{O}=\prod_{i=1}^{k_1}\{t^0_i\}\times\vk{O}_2$. \cJ{In addition, \netheo{mainthm} coincides with Theorem 7.1 in Piterbarg (1996) when $k_1=0$.}

 b)  \cJ{In the proof of \netheo{mainthm}, an extension of Pickands theorem (\cL{see} \nelem{lemsta1} below) plays an important role.
 \Hr{We remark that Pickands theorem  \bb{(see Pickands (1969))} has been rigorously proved in Piterbarg (1972).}}

% c) The definition of $\alpha(\vk{t})$-locally stationary GRF's can be extended in different ways. In terms of the behaviour of extremes, we focused, %in this paper, \cJ{on}  this extension which illustrates the ideas. Minor extensions of the \bbb{definition} lead to much complicated and difficult %to follow proofs.

Brief outline of  the paper: We give two applications of our main result in Section 2. In Section 3 we present some preliminary results.  All the proofs are relegated in Section 4 and Appendix.
%\cJ{Some technical or tedious proofs are postponed to Appendix.}

 \COM{
\BT\label{mainthm2}
Let $\{X(\tf{t}), \tf{t}\in[0,T]^{k}\}$% with
be an $\alpt$-locally stationary  GRF that satisfies the conditions of \netheo{mainthm}. Furthermore, assume  $\alpha_i(t_i)\equiv\alpha_i, i=k_1+1,\cdots,k,$ on $\tf{U}\cE{\subset}\ \R_+^{k_2}$, some compact set with positive Lebesgue measure. Then
\eqref{EM} holds with constant $\kal{K}$ given by
\BQN
\kal{K}_{\cH{\vk{U}}}&=&2^q \prod_{i=1}^{k_1}\left((C_i)^{1/\alpha_i}(B_i)^{-1/\beta_i}\right)
\cE{\prod_{i=1}^{k_1} \frac{(\alpha_i^2/2)^{1/\beta_i} \Gamma(1/\beta_i)}{\beta_i} }
\prod_{i=1}^k\mH_{\alpha_i} \left(\int_{\tf{U}}\prod_{i=k_1+1}^{k}(C_i(x_i))^{1/\alpha_i}d\tf{x}\right).
\EQN

\COM{\BQN
\pk{\underset{\tf{t}\in [0,T]^{k_1}\times\tf{U}}{\sup}X(\tf{t})>u}&=&
&&=\prod_{i=1}^{k_1}\left((C_i)^{1/\alpha_i}(B_i)^{-1/\beta_i}\right)
\mH_\alpha \left(\int_{\tf{U}}\prod_{i=k_1+1}^{k}(C_i(x_i))^{1/\alpha_i}d\tf{x}\right)
\int_{\R^{k_1}_+}e^{-\sum_{i=1}^{k_1}\frac{2}{\alpha_i^2}x_i^{\beta_i}}d\tf{x}\left(\frac{\prod_{i=1}^{k}u^{2/\alpha_i}}{\prod_{i=1}^{k_1}(\ln u)^{1/\beta_i}}\right) \Psi(u)(1+o(1))
\EQN
as $u\rw\IF$, where $C_i=C_i(0), i=\gok$.

\BQNY
&&\pk{\underset{\tf{t}\in [0,T]^{k_1}\times\tf{U}}{\sup}X(\tf{t})>u}\\
&&=\prod_{i=1}^{k_1}\left((C_i)^{1/\alpha_i}(B_i)^{-1/\beta_i}\right)\mH_\alpha \left(\int_{\tf{U}}\prod_{i=k_1+1}^{k}(C_i(x_i))^{1/\alpha_i}d\tf{x}\right)
\int_{\R^{k_1}_+}e^{-\sum_{i=1}^{k_1}\frac{2}{\alpha_i^2}x_i^{\beta_i}}d\tf{x}\left(\frac{\prod_{i=1}^{k}u^{2/\alpha_i}}{\prod_{i=1}^{k_1}(\ln u)^{1/\beta_i}}\right) \Psi(u)(1+o(1))
\EQNY
as $u\rw\IF$, where $C_i=C_i(0), i=\gok$.}
\ET
}

\section{Applications}

\cJ{In this section we apply our results to two interesting cases of $\alpha(\vk{t})$-locally stationary GRF's,
 namely, the \cL{aggregate} \bbb{mfBm's} and the $\chi$-processes generated by \cL{mfBm's} defined below.}

Let $\{B_{\alpha(t)}(t),t\ge0\}$ be a mfBm  with parameter $\alpha(t)\in(0,2], \cL{t\ge0}$. %being a H\"{o}lder function of exponent $\gamma$ such that $0<H(t)<\min(1,\gamma),t>0$.
 We define the \bbb{standardized}/normalized mfBm by
$$\overline{B}_\alpha(t)=\frac{B_{\alpha(t)}(t)}{\sqrt{Var(B_{\alpha(t)}(t))}},t\in[T_1,T_2], \quad \text{  with } 0<T_1<T_2<\IF.$$
 As shown in D\c{e}bicki and Kisowski (2008)
\BQNY
1-Cov(\ovB_\alpha(t), \ovB_\alpha(s+t))=\frac{1}{2}t^{-\alpha(t)}|s|^{\alpha(t)}+o(|s|^{\alpha(t)})
\EQNY
uniformly with respect to $t\in[T_1,T_2]$, as  $s\rw 0$.%

{\bf \underline{\cL{Aggregate} multifractional Brownian motions}}:
Let $\{\ovB_{\alpha_i}(t_i), t_i\in[T_1,T_2]\}, i=1,\cdots,k,$ be independent \bbb{standardized} mfBm's, with parameters $\alpha_i(t_i), \cL{t\ge0}, i=1,\cdots,k,$  respectively.
Assume, for any fixed $i=1,\cdots,k,$ that $\alpha_i(t_i)$ attains its minimum at the unique point $t_i^0\in (T_1,T_2)$, and that there exist some positive $M_i, \beta_i,$ and $\delta_i>1,  i=1,\cdots,k,$ \ee{such that $A2$ is satisfied. } Set $X(\tf{t})=\frac{1}{\sqrt{k}}\left(\ovB_{\alpha_1}(t_1)+\cdots+\ovB_{\alpha_k}(t_k)\right), \tf{t}\in [T_1,T_2]^k$. It  \cH{follows} that,
 as $\tf{s}\rw\tf{0}$, %, for $\tf{t}, \tf{t+s}\in [T_1,T_2]^k$,
\BQNY
1-Cov(X(\tf{t}),X(\tf{t}+\tf{s}))&=&1-\frac{1}{k}\sum_{i=1}^k Cov(\ovB_{\alpha_i}(t_i), \ovB_{\alpha_i}(t_i+s_i))\\%,\cdots,Cov(\ovB_{H_d}(t_d), \ovB_{H_d}(t_d+s_d))\right)\\
&=&\frac{1}{2k}\sum_{i=1}^k\left((t_i)^{-\alpha_i(t_i)}|s_i|^{\alpha_i(t_i)}\right)(1+o(1))
\EQNY
\def\ALB{\alpha_{i}^0}
uniformly with respect to $\tf{t}\in[T_1,T_2]^k$.  Therefore, conditions $D1-D4$ are satisfied and we have from \netheo{mainthm}
(\cJ{recall $\alpha_i:=  \alpha_i(t_i^0)$})
\BQN
&&\pk{\underset{\tf{t}\in [T_1,T_2]^k}{\sup}X(\tf{t})>u}\notag\\
&&
=2^k(2k)^{-\sum_{i=1}^k\frac{1}{ \alpha_i}} \left(\prod_{i=1}^{k}\frac{\mH_{\alpha_i}
\Gamma(1/\beta_i\ee{+1})}{t_i^0}\left(\frac{\alpha_i^2}{2M_i}\right)^{\bbb{{1/\beta_i}}}%\int_{\R^{d}_+}e^{-\sum_{i=1}^{d}\frac{ x_i^{\beta_i}}{2(H_i(t_i^0))^2}}d\tf{x}
 \right) \frac{u^{\sum_{i=1}^k\frac{2}{\alpha_i}}}{(\ln u)^{\sum_{i=1}^k 1/\beta_i}} \Psi(u)(1+o(1))
\EQN
as $u\rw\IF$.

{\bf \underline{$\chi$-processes}}:
Let $\{\ovB_{i,\alpha}(t), t\in[T_1,T_2]\}, i=1,\cdots, k,$ be independent copies of $\{\ovB_\alpha(t), t\in[T_1,T_2]\}$.
Assume that $\alpha(t)$ attains its  minimum at the unique point $t^0\in (T_1,T_2)$, and that there exist some positive $M, \beta,$ and $\delta>1,$ such that again $A2$ holds.  Consider the $\chi$-process defined by
$$\chi_k(t)=\sqrt{\ovB_{1,\alpha}^2(t)+ \cdots + \ovB_{k,\alpha}^2(t)}, \quad t\in[T_1,T_2].$$
Further, we introduce a  GRF $$
Y(t,\tf{u})=\ovB_{1,\alpha}(t)u_1+\cdots+\ovB_{k,\alpha}(t)u_k, \quad {\tf{u}=(u_1,\cdots,u_k)}
$$
defined on the cylinder $\mathcal{G}_T=[T_1,T_2]\times\mathcal{S}_{k-1}$, with $\mathcal{S}_{k-1}$ being the unit sphere in $\R^k$ (with respect to $L_2$-norm). In the light of Piterbarg (1996)
$$
\underset{t\in[T_1,T_2]}\sup \chi_k(t)=\underset{(t,\tf{u})\in\mathcal{G}_T}\sup Y(t,\tf{u}).
$$
%We can
\Hr{Further we have as $(s,\tf{v})\rw\cL{(0,\tf{0})}$}
\BQNY
1-Cov(Y(t,\tf{u}),Y(t+s,\tf{u+v}))
%&=&1-\sum_{i=1}^k Cov(\ovB_{\alpha}(t),\ovB_\alpha(t+s))u_i(u_i+v_i)\\
%&=&1-Cov(\ovB_{\alpha}(t),\ovB_\alpha(t+s))%\sum_{i=1}^k\left(\frac{(u_i)^2+(u_i+v_i)^2}{2}\right)\\
%+\frac{1}{2}Cov(\ovB_{\alpha}(t),\ovB_\alpha(t+s))\sum_{i=1}^k\abs{v_i}^2\\
&=&\frac{1}{2}t^{-\alpha(t)}|s|^{\alpha(t)}+\frac{1}{2}\sum_{i=1}^{\cL{k-1}} \abs{v_i}^2+o\left(|s|^{\alpha(t)}+\sum_{i=1}^{\cL{k-1}} \abs{v_i}^2\right)
\EQNY
uniformly with respect to $(t,\tf{u})\in \mathcal{G}_T$.
Therefore, the conditions $D1-D4$ are satisfied and we have that \cH{(recall Remark a) above)}
\BQN
\pk{\underset{{t}\in [T_1,T_2]}{\sup}\chi_k({t})>u}= 2^{\cL{\frac{5}{2}}-\frac{k}{2}-\frac{1}{\beta}-\frac{1}{\alpha(\cH{t^0})}} \frac{\mH_{\alpha(t^0)} (\alpha(\cH{t^0}))^{\frac{2}{\beta}}\Gamma(\frac{1}{\beta}\cJ{+1})}{M^{1/\beta}\cH{t^0} \Gamma({\cH{k}/2})}
\frac{u^{\cL{k-1}+\frac{2}{\alpha(t^0)}}}{(\ln u)^{1/\beta}} \Psi(u)(1+o(1)), \ u\rw\IF.
\EQN

\section{\cL{Preliminary Lemmas}}
\cL{This section is concerned with some preliminary lemmas used for the proof of \netheo{mainthm}.}
%We assume, {without loss}, that $1\le k_1<k$, $M_i=1$ for any $i=1,\cdots, k_1$, and $t^0_i={0}$ for any $i=1,\cdots,k_1$.
{We assume, {without loss}, that $1\le k_1<k$ and $M_i=1, i=1,\cdots, k_1$.  As pointed out in D\c{e}bicki and Kisowski (2008), for the asymptotics of the original process, we have to replace
$C_i(\cdot)$ with $(M_i)^{-\alpha_i/\beta_i} C_i(\cdot),\  i=1,\cdots, k_1.$ We may further assume that $t^0_i={0}, i=1,\cdots,k_1$, and thus the %for the general cases, the
final general result should be multiplied by $2^{q_{k_1}}$.}
\cL{Hereafter, consider $\{X(\tf{t}), \tf{t}\in[0,T]^k\}$ to be an $\alpt$-locally stationary GRF with the above simplification (called {\it simplified $\alpt$-locally stationary GRF}).} % that satisfies the  conditions $A1$-$A2$  with the assumptions above and \eqref{eqcov11}.
Set next
\BQNY
t_u^i= \left(\frac{(\alpha_i)^2}{\beta_i}\frac{\ln\ln u}{\ln u}\right)^{\frac{1}{\beta_i}}, \quad i=1,\cdots,k_1.
\EQNY
Clearly
%It  can be  shown that (see \nelem{lemtail} \cH{below})
\BQN\label{equplow}
&&\pk{\underset{\tf{t}\in \prod_{i=1}^{k_1}[0,t^i_u] \times\prod_{i=k_1+1}^k[a_i,b_i]}{\sup}X(\tf{t})>u}\le\pk{\underset{\tf{t}\in [0,T]^k}{\sup}X(\tf{t})>u}\le\nonumber\\
&&\ \ \  \le\pk{\underset{\tf{t}\in \prod_{i=1}^{k_1}[0,t^i_u] \times\prod_{i=k_1+1}^k[a_i,b_i]}{\sup}X(\tf{t})>u}
 +\pk{\underset{\tf{t}\in \left([0,T]^k/\prod_{i=1}^{k_1}[0,t^i_u] \times\prod_{i=k_1+1}^k[a_i,b_i]\right)}{\sup}X(\tf{t})>u}.
 %(1+o(1)), %\quad \mathrm{as}\quad u\rightarrow \infty,
\EQN
There are two steps in the proof of \netheo{mainthm}. In step 1, we focus on the asymptotics of
\BQN
\bb{\Pi(u)}:=\pk{\underset{\tf{t}\in \prod_{i=1}^{k_1}[0,t^i_u] \times\prod_{i=k_1+1}^k[a_i,b_i]}{\sup}X(\tf{t})>u}, \quad u\rightarrow \infty,
\label{localsta2}
\EQN
which is the main part of our proof. In step 2, we shall show that (see \nelem{lemtail}  below)
\BQN\label{step2}
\pk{\underset{\tf{t}\in \left([0,T]^k/\prod_{i=1}^{k_1}[0,t^i_u] \times\prod_{i=k_1+1}^k[a_i,b_i]\right)}{\sup}X(\tf{t})>u}=o(\Pi(u)), \ \ \ %\quad \mathrm{as}\quad
u\rightarrow \infty.
\EQN
\cL{The idea of finding the asymptotics of \eqref{localsta2} is based on the so-called  double-sum method; see e.g., Pickands (1969) or Piterbarg (1996). Before going to the detail of the proof, let us recall the brief outline of the double-sum method. First of all, we need to find a suitable partition, say cubes $\{W_u^i\}$, of the set $\prod_{i=1}^{k_1}[0,t^i_u] \times\prod_{i=k_1+1}^k[a_i,b_i]$. Then using the well-known Bonferroni's inequality we find upper and lower bounds for \eqref{localsta2}, i.e.,}
\BQNY
&&\sum_i\pk{\underset{\tf{t}\in W_u^i}{\sup}X(\tf{t})>u}\ge\pk{\underset{\tf{t}\in \prod_{i=1}^{k_1}[0,t^i_u] \times\prod_{i=k_1+1}^k[a_i,b_i] }{\sup}X(\tf{t})>u}\ge\\
 &&\ \ \ \ \ \ \ \ \ \ \  \ge\sum_i\pk{\underset{\tf{t}\in W_u^i}{\sup}X(\tf{t})>u}- \underset{i<j}{\sum\sum}\pk{\underset{\tf{t}\in W_u^i}{\sup}X(\tf{t})>u, \underset{\tf{t}\in W_u^j}{\sup}X(\tf{t})>u}.
\EQNY
\cL{Finally, we show that the asymptotics of the single-sum terms on both sides are the same and the double-sum term is relatively negligible. In what follows, we shall first introduce the cubes that \Hr{are} used as the partition, followed then by some preliminary results (Lemmas \ref{lemslep1}-\ref{lemdoub}) concerning the estimation for the summands of both single-sum and double-sum terms in the last formula.}
\bb{For $i=1,\cdots,k_1,$ set}
$$c_{p_i}^i=c_{p_i}^i(u):=\left(\frac{p_i}{\ln u(\ln \ln u)^{1/\beta_i}}\right)^{1/\beta_i}, \quad
 A_{p_i}^i=A_{p_i}^i(u):=[c_{p_i}^i,c_{p_i+1}^i],$$
 and let $m_i=m_i(u):=\lfloor\frac{(\alpha_i)^2}{\beta_i}(\ln\ln u)^{1+1/\beta_i}\rfloor$, where $\lfloor x \rfloor$ denotes the integer part of $x$. Further, let $S>1$ \bb{be a fixed constant;} by dividing each $A_{p_i}^i$ into subintervals of length $S/u^{2/(\alpha_i(c_{p_i+1}^i))}$ (recall function $\alpha_i(\cdot)$ in \eqref{eqcov1}), we define
\BQNY
B_{{j_i},p_i}^i=B_{{j_i},p_i}^i(u):=\left[c_{p_i}^i+\frac{{j_i} S}{u^{2/(\alpha_i(c_{p_i+1}^i))}}, c_{p_i}^i+\frac{({j_i}+1) S}{u^{2/(\alpha_i(c_{{p_i}+1}^i))}}\right]
\EQNY
for ${j_i}=0,1,\cdots,n_{i,{p_i}}=n_{i,{p_i}}(u):=\lfloor \frac{c^i_{{p_i}+1}-c^i_{p_i}}{S}u^{2/(\alpha_i(c_{{p_i}+1}^i))}\rfloor$.

 Moreover, let $k_2:=k-k_1$, $\tf{a}=(a_{k_1+1},\cdots,a_k),$ \cE{and let} $\kk=(K_1,\cdots,K_{k_2})\in \mathbb{Z}^{k_2}$ be a vector with integer coordinates. %\cH{We assume for simplicity in the following that $k-k_1\ge 1$.}
 For $\delta>0$, we denote
  $$\delta_{\kk}=(\tf{a}+\delta\kk+[0,\delta]^k)\cap \prod_{i=k_1+1}^k[a_i,b_i],$$
   where $\kk\in \mathcal{B}$ with
\BQNY
\mathcal{B}&=&\{\kk\in\ZZ^{k_2}:  \delta_{\kk}\neq \emptyset \}.
\EQNY
Define an operator $g_u$  on $\R^{k_2}$ as in Piterbarg (1996), i.e., for $\tf{t}=(t_{k_1+1},\cdots,t_{k})\in\R^{k_2}$
\BQN
g_u\tf{t}=\left(u^{-\frac{2}{\alpha_{k_1+1}}}t_{k_1+1},\cdots,u^{-\frac{2}{\alpha_k}}t_k\right).\label{eqgu}
\EQN
Denote $\triangle_0=g_u[0,1]^{k_2},$ and, for fixed $\kk\in\mathcal{B}$, $\triangle_{\ii}=\triangle_{\ii}(u):=g_uS\ii +\triangle_0S$ with $\ii=(I_1^{\kk},\cdots,I_{k_2}^{\kk})\in\ZZ^{k_2}$ being a vector with integer coordinates. Further, let $V_{\ii,\kk}:=\tf{a}+\delta\kk+\triangle_{\ii}$, where $\ii\in \mathcal{A}_{\kk}$ with
\BQNY
\mathcal{A}_{\kk}=\{\ii\in\ZZ^{k_2}:   V_{\ii,\kk}\cap \delta_{\kk}\neq \emptyset\}.
\EQNY
Denote
$$N_{\kk}^+=\# \{\ii\in\ZZ^{k_2}:   V_{\ii,\kk}\cap \delta_{\kk}\neq \emptyset\}\ \mq{and}\  N_i=\left\lfloor\frac{\delta}{S}u^{2/\alpha_{k_1+i}}\right\rfloor, i=\gokkk.$$  %$N_{\kk}^-=\# \{\ii:   V_{\ii,\kk}\subset \delta_{\kk}\}$.
Moreover, let, for $i=1,\cdots,k_1$,
\BQNY
\LL{1}^i&=&\{({j_i},{p_i}): {j_i},{p_i}\in\mathbb{Z}, 0\le {p_i}\le m_i-1, 0\le {j_i}\le n_{i,{p_i}}-1 \},\\
\UU{1}^i&=&\{({j_i},{p_i}):  {j_i},{p_i}\in\mathbb{Z}, 0\le {p_i}\le m_i, 0\le {j_i}\le n_{i,{p_i}} \},
\EQNY
and
\BQNY
\LL{2}=\{(\ii,\kk): \kk\in \BB, V_{\ii,\kk}\subset \delta_{\kk} \},\  \ \
\UU{2}=\{(\ii,\kk):  \kk\in \BB, \ii\in\AAk\}.
\EQNY
\bb{We} have
\BQNY
\underset{(\ii,\kk)\in  \LL{2}}{\underset{(j_i,p_i)\in\LL{1}^i,i=1,\cdots,k_1 }\bigcup} \prod_{i=1}^{k_1}B_{j_i,p_i}^i\times V_{\ii,\kk}\ \subset \  \prod_{i=1}^{k_1}[0,t^i_u] \times\prod_{i=k_1+1}^k[a_i,b_i]\ \subset\
\underset{ (\ii,\kk)\in  \UU{2}}{\underset{(j_i,p_i)\in\UU{1}^i,i=1,\cdots,k_1}\bigcup} \prod_{i=1}^{k_1}B_{j_i,p_i}^i\times V_{\ii,\kk}.
\EQNY
In order to specify the 'distance' between segments of the type $\prod_{i=1}^{k_1}B_{j_i,p_i}^i\times V_{\ii,\kk}$, we introduce the following order relation: for any $(j,p), (j',p')\in \ZZ^2$, we write
\BQNY
%&&(j,p)\preceq (j',p')\quad \mathrm{iff}\quad (p< p')\ \mathrm{or}\ (p=p'\ \mathrm{and}\ j\le j')\\
&&(j,p)\prec (j',p')\quad \mathrm{iff}\quad (p<p')\ \mathrm{or}\ (p=p'\ \mathrm{and}\ j<j').
\EQNY
Further, for $\tf{j},\tf{p}, \tf{j}',\tf{p}' \in \ZZ^{k_1}$ with $(j_i,p_i),(j'_i,p'_i)\in \LL{1}^i, i=1,\cdots,k_1,$
\BQNY
%&&(\tf{j},\tf{p})\preceq (\tf{j}',\tf{p}') \quad \mathrm{iff}\quad (j_i,p_i)\preceq (j'_i,p'_i) \quad \mathrm{for\ all}\quad i=1,\cdots,k_1,\\
&&(\tf{j},\tf{p})\prec (\tf{j}',\tf{p}') \quad \mathrm{iff}\quad (j_i,p_i)\prec (j'_i,p'_i) \quad \mathrm{for\ some}\ i=1,\cdots,k_1,
\mathrm{and}\ (j_l,p_l)= (j'_l,p'_l)\ \mathrm{for}\ l=1,\cdots,i-1,
\EQNY
and, for $(\ii,\kk), (\iii,\kk')\in \LL{2},$%(\mathrm{or}\ \UU{2})$,
\BQNY
%&&(\ii,\kk)\preceq (\ii',\kk') \quad \mathrm{iff}\quad (I^{\kk}_i,K_i)\preceq (I'^{\kk}_i, K'_i) \quad \mathrm{for\ all}\quad i=1,\cdots,k_2\\
&&(\ii,\kk)\prec (\iii,\kk') \quad \mathrm{iff}\ (I^{\kk}_i,K_i)\prec (I'^{\kk'}_i, K'_i) \  \mathrm{for \ some}\ i=1,\cdots,k_2, \mathrm{and}\
 (I^{\kk}_l,K_l)=(I'^{\kk'}_l, K'_l) \ \mathrm{for}\ l=1,\cdots,i-1.
\EQNY
Moreover, define, for ${j,p},{j',p'}\in \mathbb{Z}$,
$$
N_{j,p}^{j',p'}:=\#\{(j'',p'')\in \mathbb{Z}^2: (j,p)\prec(j'',p'')\prec(j',p')\}.
$$
\cJ{In the sequel}, for  fixed $j_i,p_i, \ii, \kk$ \cJ{such that $(j_i,p_i)\in \UU{1}^i, i=1,2,\cdots,k_1$ and  $(\ii,\kk)\in\UU{2}$,} we consider the  GRF $X(\tf{v}):=X(v_1,\cdots,v_{k})$ on
 $$\Ajpik:=\prod_{i=1}^{k_1}B_{j_i,p_i}^i\times V_{\ii,\kk}.$$
%\cE{It is better to prove 3.8 and 3.9 in a Lemma (we can put in Appendix, for instance).}
%\cE{Since Lemma 3.2 is very technical and special for our case, we might integrate it in the proof and "delete" Lemma 3.2 but state the results in the proof of Lemma 3.4 below.}
\cL{In order to obtain the estimates of the tail probabilities of the supremum of $X$ on $\Ajpik$ (see Lemmas \ref{lemslep1} and \ref{lemsta2} below),} we introduce the following \cL{stationary} GRF's, for a fixed (marked) point $\tf{v}^0=(v^0_1,\cdots,v_k^0):=\tf{v}_{\vk{j,p},\ii,\kk}^0$ in $\Ajpik$ :\\
----$\{Y_{\vn,u}^{\tf{v}^0}(\tf{\nu}), \tf{\nu}\in[0,S]^{k}\}$ is  a family of centered stationary  GRF's   with
\BQNY
Cov(Y_{\vn,u}^{\tf{v}^0}(\tf{\nu}),Y_{\vn,u}^{\tf{v}^0}(\tf{\nu}+\tf{x}))\ =\ e^{-(1-\vn)\left(\sum_{i=1}^{k_1} C_i(\tf{v}^0)u^{-2}|x_i|^{\alpha_i+2(t_u^i)^{\beta_i}}
+\sum_{i=k_1+1}^{k} C_i(\tf{v}^0)u^{-2}|x_i|^{\alpha_i}\right)}
\EQNY
for $\vn\in(0,1)$, $u>0$ such that $\alpha_i+2(t_u^i)^{\beta_i}\le2$, $i=1,\cdots,k_1$, and $\tf{\nu,\nu}+\tf{x}\in[0,S]^{k}$.

----$\{Z_{\vn,u}^{\tf{v}^0}(\tf{\nu}), \tf{\nu}\in[0,S]^{k}\}$ \cH{is} a family of centered stationary  GRF's   with
\BQN
Cov(Z_{\vn,u}^{\tf{v}^0}(\tf{\nu}),Z_{\vn,u}^{\tf{v}^0}(\tf{\nu}+\tf{x}))\ =\ e^{-(1+\vn)\left(%\sum_{i=1}^{k_1} C_i(\tf{v}^0)u^{-2}(x_i)^{\alpha_i}
\sum_{i=1}^{k} C_i(\tf{v}^0)u^{-2}|x_i|^{\alpha_i}\right)},\label{eqcovZ}
\EQN
for $\vn>0$, $u>0$  and $\tf{\nu, \nu}+\tf{x}\in[0,S]^{k}$.

\begin{lem}
\label{lemslep1}
For any $\vn\in(0,1)$, there exists $u_\vn>0$ such that for $u>u_\vn$,
\BQN
&(i)&\ \pk{\underset{\tf{v}\in \Ajpik}\sup X(\tf{v})>u}\ \ge\ \pk{\underset{\tf{\nu}\in[0,S]^{k}}\sup Y_{\vn,u}^{\tf{v}^0}(\tf{\nu})>u},\nonumber\\
&(ii)&\ \pk{\underset{\tf{v}\in \Ajpik}\sup X(\tf{v})>u}\ \le\ \pk{\underset{\tf{\nu}\in[0,S]^{k}}\sup Z_{\vn,u}^{\tf{v}^0}(\tf{\nu})>u}.\label{slep2}
\EQN
\end{lem}
%The proof of \nelem{lemslep1} is  \bb{very} technical, and \cH{therefore is} relegated to Appendix.
\BRM
 Due to  continuity of the functions $C_i(\cdot),i=1,\cdots,k$, the point $\tf{v}^0$ can also be chosen as a fixed (marked) point in $\prod_{i=1}^{k_1}A_{p_i}^i\times \delta_{\kk}$ when $\delta$ is sufficiently small and $u$ is sufficiently large. In the sequel, we chose $\tf{v}^0$ in this way. Actually $\tf{v}^0$ depends on $\tf{p}, \kk$, but, if no confusion is caused, \cH{for} notational simplicity
 we still \cH{write} $\tf{v}^0$.
\ERM
\bbb{Next we introduce a {\it structural modulus} on $\R^k$ by
 %for, we set %introduce the following notation:
$$
{|\tf{s}|}_{\vk{\alpha}}=\sum_{i=1}^k{|s_i|}^{\alpha_i}, \ \ \tf{s}\in\R^k. \ %\vk{\alpha}=(\alpha_1,\cdots,\alpha_k).
$$}
The following \bb{result inspired} by Lemma 7 of H\"{u}sler and Piterbarg (2004) is crucial for our investigation;
\bb{its proof is  relegated to Appendix.}
%, which is important for us to get the exact asymptotics of the  GRF's   considered.
%This is an independent interesting result, but since the proof of it is tedious, it is relegated to Appendix.
\begin{lem}\label{lemsta1}
For any compact set $\tf{D}\in \R_+^{k}$, let $\{X_u(\tf{t}), \tf{t}\in \tf{D}\}$, $u>0,$ be a family of a.s.\ continuous  GRF's  , with $\E{X_u(\tf{t})}\equiv0$, $\E{(X_u(\tf{t}))^2}\equiv1$ for all $u$, and with correlation function $r_u(\tf{t},\tf{s})=\E{X_u(\tf{t})X_u(\tf{s})}.$ If
\BQN
\underset{u\rw\infty}\lim u^2(1-r_u(\tf{t},\tf{s}))={|\tf{t}-\tf{s}|}_{\tf{\alpha}}\label{eqlimru}
\EQN
uniformly with respect to $\tf{t}, \tf{s}\in \tf{D}$,
%(2) For  some structure $(L,\beta)$ and $C$ and all $\tf{t}, \tf{s}\in\tf{D}, u\ge0$, $u^2(1-r_u(g_u\tf{t},g_u\tf{s}))\le C|\tf{t}-\tf{s}|_{L,\beta}$,
then %for any $\tf{D}\in\R_+^d$,
\BQNY
\pk{\underset{\tf{t}\in \tf{D}}\sup X_u(\tf{t})>u}= \mH_{(k, \tf{\alpha})}[\tf{D}]\Psi(u)(1+o(1))
\EQNY
as $u\rw\infty$, where
\BQN
\mH_{(k, \tf{\alpha})}[\tf{D}]=\E{\exp\left(\underset{\tf{t}\in\tf{D}}\sup (\widetilde{B}_{\tf{\alpha}} (\tf{t})-|\tf{t}|_{\tf{\alpha}})\right)}\in(0,\IF)
\label{PKD}
\EQN
as defined in Piterbarg (1996), with
$$
\bb{\widetilde{B}}_{\tf{\alpha}}(\tf{t})=\bbb{\sqrt{2}}\sum_{i=1}^k B_{\alpha_i}^{(i)}(t_i)
$$
and $B_{\alpha_i}^{(i)}, 1\le i\le k$, being independent \bb{fBm's} %fractional Brownian motions
with Hurst \Hr{indexes} $\alpha_i/2 \cL{\in(0,2]}$, respectively.
\end{lem}

\begin{lem}\label{lemsta2}
For any $S>1$ and $\vn\in(0,1)$, we have, as $u\rw\infty,$

$(i)\ \  \pk{\underset{\tf{\nu}\in[0,S]^{k}}\sup Y_{\vn,u}^{\tf{v}^0}(\tf{\nu})>u}=\prod_{i=1}^{k}\mH_{\alpha_i}\left[0,(\cJ{C_i(\tf{v}^0)}(1-\vn))^{1/\alpha_i}S\right]\Psi(u)(1+o(1)),$

$(ii)\ \  \pk{\underset{\tf{\nu}\in[0,S]^{k}}\sup Z_{\vn,u}^{\tf{v}^0}(\tf{\nu})>u}=\prod_{i=1}^{k}\mH_{\alpha_i}\left[0,(\cJ{C_i(\tf{v}^0)}(1+\vn))^{1/\alpha_i}S\right]\Psi(u)(1+o(1))$,

where \cJ{(recall \eqref{PKD}) %, for any positive number $S,$
%\BQNY
\Hr{we set} $\mH_{\alpha_i}[0,S]:=\mH_{(1,\alpha_i)}[[0,S]],\ \  i=1,2,\cdots,k.$
}%\EQNY}
\end{lem}

In order to estimate the double-sum term in the derivation of \eqref{localsta2}, we \Hr{need} the following two lemmas.

\BEL\label{lemdoubZ}
Let \bb{GRF} $\{\widetilde{Z}^{\tf{w}^0}_{\tf{\vn},u}(\nu); \tf{\nu}\in[0,S]^k\}$, having  covariance structure  \eqref{eqcovZ} with $\tf{v}^0$ replaced by $\tf{w}^0$, \cE{be independent} of $\{Z_{\vn,u}^{\tf{v}^0}(\tf{\nu}); \tf{\nu}\in[0,S]^k\}$, with $\vn>0$. Then there exists some positive constant ${F}_\vn$, \cJ{for $u$ large enough,} we have
\BQNY
\pk{\underset{\tf{\nu,\mu}\in[0,S]^{k}}\sup \frac{1}{\sqrt{2}}\left(Z_{\vn,u}^{\tf{v}^0}(\tf{\nu})+\widetilde{Z}^{\tf{w}^0}_{\vn,u}(\tf{\mu})\right)>u}\le F_\vn S^{2k}\Psi(u).%(1+o(1))
%&&=\left(\prod_{i=1}^{k}\mH_{\alpha_i}\left[0,(C_i(\tf{v}^0) (1+\vn))^{1/\alpha_i}S\right]\right)\left(\prod_{i=1}^{k}\mH_{\alpha_i}\left[0,(C_i(w^0_i) %(1+\vn))^{1/\alpha_i}S\right]\right)\Psi(u)(1+o(1)).
\EQNY
%where $C_i=C_i(0), i=1,\cdots,k$.
\EEL

%In order to estimate the double-sum term in the proof of \netheo{mainthm}, we give the following lemma.
Next, we introduce a distance of two sets $\tf{D}_1, \tf{D}_2\subset \cH{\R}_+^k$ by
$$
dist(\tf{D}_1,\tf{D}_2)=\underset{\tf{t}\in\tf{D}_1, \tf{s}\in\tf{D}_2}\inf|\tf{t}-\tf{s}|_{\tf{\alpha}}.
$$
Further, we fix some sufficiently small $\gamma_0>0$ in the following way:
uniformly with respect to $\tf{t}\in [0,T]^k$,
\BQN
1-Cov(X(\tf{t}),X(\tf{t}+\tf{s}))<\eta_0\in[0,1/2)\label{eqcovvn}
\EQN
for $|\tf{s}|_{\tf{\alpha}}<\gamma_0$ (recall \eqref{eqcov1}).

\BEL\label{lemdoub}
There exist some universal positive constants $\mathbb{C}, \mathbb{C}_1$ such that, \cJ{for sufficiently large $u,$} the following statements are established.

$(1)$ For $(j_i,p_i),(j'_i,p'_i)\in \LL{1}^i, i=1,\cdots,k_1,(\ii,\kk),(\iii,\kk')\in\LL{2}$ satisfying
\BQN
dist\left(\Ajpik,\Ajpikp\right)< \gamma_0\label{eqdis}
\EQN
 and
\COM{
\qquad $(a)$ $(\tf{j}, \tf{p})\prec (\tf{j}', \tf{p}')$ and
$$N_{j_i,p_i}^{j'_i,p'_i}>0 \ \mathrm{for\ some\ }i=1,\cdots,k_1, \ \mathrm{or}\ N_{I^{\kk}_i,K_i}^{I'^{\kk}_i,K'_i}>0 \ \mathrm{for\ some\ }i=1,\cdots,k_2,$$

\qquad $(b)$ $(\ii,\kk)\preceq(\iii,\kk')$,
}
$$N_{j_i,p_i}^{j'_i,p'_i}>0 \ \mathrm{for\ some\ }i=1,\cdots,k_1, \ \mathrm{or}\ N_{I^{\kk}_i,K_i}^{I'^{\kk'}_i,K'_i}>0 \ \mathrm{for\ some\ }i=1,\cdots,k_2,$$
 we have %as $u\rw\IF$,
\BQN
&&\pk{\underset{\tf{v}\in \Ajpik }\sup X(\tf{v})>u, \underset{\tf{v'}\in \Ajpikp }\sup X(\tf{v'})>u}\notag\\
&&\le \mathbb{C}S^{2k}\expo{{-\mathbb{C}_1\left(\sum_{i=1}^{k_1}\left(\sqrt{N_{j_i,p_i}^{j'_i,p'_i}}S\right)^{\alpha_i}+\sum_{i=1}^{k_2}
\left(N_{I^{\kk}_i,K_i}^{I'^{\kk'}_i,K'_i}S\right)^{\alpha_{k_1+i}}\right)}}\Psi(u).\label{eqdoub1}
\EQN
%where
%\cE{$$ A=  A_{\vk{j,p},\ii,\kk}, \quad  A'=\prod_{i=1}^{k_1}B_{j'_i,p'_i}^i\times V_{\ii',\kk'}.$$}

$(2)$ \cE{Let}  $(j_i,p_i),(j'_i,p'_i)\in \LL{1}^i, i=1,\cdots,k_1,(\ii,\kk),(\iii,\kk')\in\LL{2}$ satisfy %\eqref{eqdis} and
%$(\tf{j}, \tf{p})\preceq(\tf{j}', \tf{p}'),(\ii,\kk)\preceq(\ii',\kk')$ and
%$(\tf{j}, \tf{p})\prec(\tf{j}', \tf{p}')\ \mathrm{or}\ (\ii,\kk)\prec(\ii',\kk')$ and
%and
$$N_{j_i,p_i}^{j'_i,p'_i}=0 \ \mathrm{for\ all\ }i=1,\cdots,k_1, \ \mathrm{and}\ N_{I^{\kk}_i,K_i}^{I'^{\kk'}_i,K'_i}=0 \ \mathrm{for\ all\ }i=1,\cdots,k_2.$$

If $(\tf{j}, \tf{p})\prec(\tf{j}', \tf{p}')$, then the following number $\kappa$ can be defined:
\BQNY
\kappa=\left\{
              \begin{array}{ll}
i_1^1:=\inf\{1\le i \le k_1: p_i=p'_i, j'_i=j_i+1\}, &  \mathrm{if}\ i_1^1\exists,\\
i_2^1:=\inf\{1\le i \le k_1: p'_i=p_i+1, j_i=n_{i,p_i},j'_i=0\}, &  \ \mathrm{if}\ i_1^1\not\exists.% i_2^1\exists,
              \end{array}
            \right.
\EQNY
Similarly, if $(\tf{j}, \tf{p})=(\tf{j}', \tf{p}')$ and $(\ii,\kk)\prec(\ii',\kk')$, then we can define $\kappa$ as
\BQNY
\kappa=\left\{
              \begin{array}{ll}
i_1^2:=k_1+\inf\{1\le i\le k_2:K_i=K'_i,I'^{\kk'}_i=I^{\kk}_i+1\}, &  \ \mathrm{if}\ i_1^2\exists,\\
i_2^2:=k_1+\inf\{1\le i\le k_2:K'_i=K_i+1, I^{\kk}_i=N_{i},I'^{\kk'}_i=0\}, &  \ \mathrm{if}\ i_1^2\not\exists.
              \end{array}
            \right.
\EQNY
Assume, without loss of generality, that $\kappa=i_1^1$ exists. We have%, as $u\to \IF$,
\BQN
%&&
\pk{\underset{\tf{v}\in \Ajpik }\sup X(\tf{v})>u, \underset{\tf{v'}\in A''_\kappa }\sup X(\tf{v'})>u}%\\
&\le& \mathbb{C}S^{2k}\expo{{-\mathbb{C}_1S^{\alpha_{\kappa}/2}}}\Psi(u),\label{eqdoub2}
\EQN
where %as $u\rw\IF$.
$$ A''_\kappa = \prod_{i=1}^{\kappa-1}B_{j'_i,p'_i}^i\times \left[c_{p_\kappa}^\kappa+\frac{{(j_\kappa+1)} S+\sqrt{S}}{u^{2/(\alpha_\kappa(c_{p_\kappa+1}^\kappa))}}, c_{p_\kappa}^\kappa+\frac{({j_\kappa}+2) S}{u^{2/(\alpha_\kappa(c_{{p_\kappa}+1}^\kappa))}}\right]   \times\prod_{i=\kappa+1}^{k_1}B_{j'_i,p'_i}^i\times V_{\ii',\kk'}.$$

$(3)$ If $(j_i,p_i),(j'_i,p'_i)\in \LL{1}^i, i=1,\cdots,k_1,(\ii,\kk),(\iii,\kk')\in\LL{2}$ satisfy
\BQN
dist\left(\Ajpik,\Ajpikp\right)\ge \gamma_0,\label{eqdis2}
\EQN
 then there exist some \bbb{constants (independent of $u$)} $h>0$ and $\lambda\in(0,1)$ such that
\BQN
\pk{\underset{\tf{v}\in \Ajpik }\sup X(\tf{v})>u, \underset{\tf{v'}\in \Ajpikp }\sup X(\tf{v'})>u}\le 2\Psi\left(\frac{u-h/2}{\sqrt{1-\lambda/2}}\right).\label{eqdoub3}
\EQN
%for sufficiently large $u$.
\EEL

\cL{The next lemma gives the asymptotics of \eqref{localsta2}, which is the main part of the proof of \netheo{mainthm}.}

\BEL\label{lemtail1} Let $\{X(\tf{t}), \tf{t}\in[0,T]^k\}$ be the simplified $\alpt$-locally stationary GRF. We have
\BQNY
\Pi(u)&=& \biggl(\prod_{i=1}^{k_1}
\Bigl( \frac{\alpha_i^2}{2 }\Bigr)^{1/\beta_i}  \cH{\Gamma(1/\beta_i+1)} \biggr)
\Bigr(\prod_{i=1}^k\mH_{\alpha_i} \Bigl)
\int_{\vk{x}\in \prod_{i=1}^{k_1}\{t^0_i\}\times\prod_{i=k_1+1}^k [a_i,b_i]} \cJ{\prod_{i=1}^{k}(C_i(\tf{x}))^{1/\alpha_i}d\tf{x}}\\
&&\times \cL{ u^{\alpha}(\ln u)^\beta} \Psi(u)(1+o(1)), \ \quad u\rw\IF,
\EQNY
 where $\alpha,\beta$ are the same as in \netheo{mainthm}.
\EEL
\def\OE{\bb{\overline{\ve}}}

The last lemma stated below establishes Eq. \eqref{step2}.

\BEL\label{lemtail}
Let $\{X(\tf{t}), \tf{t}\in[0,T]^k\}$ be the simplified $\alpt$-locally stationary GRF. Then %If $A1$ and $A2$ are satisfied, then, as $u\rw\IF$,
\BQNY
\pk{\underset{\tf{t}\in \left([0,T]^k/\prod_{i=1}^{k_1}[0,t^i_u] \times\prod_{i=k_1+1}^k[a_i,b_i]\right)}{\sup}X(\tf{t})>u}=o(\Pi(u)), \ \ u\rw\IF.
\EQNY
\EEL

\section{Proofs}

\prooftheo{mainthm}  Taking into account of the (simplification) statement in the beginning of Section 3, we conclude that the claim follows directly from \eqref{equplow} and Lemmas \ref{lemtail1} and \ref{lemtail}. \QED

\prooflem{lemslep1}
Set
$$\Xjpik(\tf{\nu})=X\left(c_{p_1}^1+\frac{j_1 S+\nu_1}{u^{2/(\alpha_1(c_{p_1+1}^1))}},\cdots,c_{p_{k_1}}^{k_1}+\frac{j_1 S+\nu_{k_1}}{u^{2/(\alpha_{k_1}(c_{p_{k_1}+1}^{k_1}))}},
\tf{a}+\delta\kk+g_uS\ii +\triangle_0^{\tf{\nu}}\right),$$
 with $\triangle_0^{\tf{\nu}}=g_u\prod_{i=k_1+1}^{k}[0,\nu_i]$. \cE{It follows that}
\BQN
\underset{\tf{v}\in  \Ajpik}\sup X(\tf{v})\ \overset{d}=\  \underset{\tf{\nu}\in[0,S]^{k}}\sup \Xjpik(\tf{\nu}).
\EQN
 Furthermore, we derive, for the fixed point $\tf{v}^0$ in $\Ajpik$, and $u$ sufficiently large,
\BQNY
&&1-Cov\left(\Xjpik(\tf{\nu}),\Xjpik(\tf{\nu}+\tf{x})\right)\\
&\ge&(1-\vn/4)^{1/3}\left(\sum_{i=1}^{k_1} C_i(\tf{v})|u^{-2/(\alpha_i(c_{p_i+1}^i))} x_i|^{\alpha_i\left(c_{p_i}^i+\frac{j_i S+\nu_i}{u^{2/(\alpha_i(c_{p_i+1}^i))}}\right)}
+\sum_{i=k_1+1}^{k} C_i(\tf{v})u^{-2}|x_i|^{\alpha_i}\right)\\
&\ge&(1-\varepsilon/2)^{1/3}\left(\sum_{i=1}^{k_1} C_i(\tf{v}^0)|u^{-2/(\alpha_i(c_{p_i+1}^i))} x_i|^{\alpha_i\left(c_{p_i}^i+\frac{j_i S+\nu_i}{u^{2/(\alpha_i(c_{p_i+1}^i))}}\right)}
+\sum_{i=k_1+1}^{k} C_i(\tf{v}^0)u^{-2}|x_i|^{\alpha_i}\right)
\EQNY
uniformly with respect to $\tf{\nu, \nu+x}\in[0,S]^k$, where we used the fact that $C_i(\cdot), i=1,\cdots,k$, are continuous functions.%as $\tf{x}\rw \tf{0}$ uniformly with respect to $\nu\in[0,S]^{k_1}\times[0,U]^{k_2}$.

\Hr{In view of the proof} of Lemma 4.1 of D\c{e}bicki and Kisowski (2008) for sufficiently large $u$ \Hr{we obtain} %\cE{and some $\ve\in (0,1)$}
\BQN
&&1-Cov\left(\Xjpik(\tf{\nu}),\Xjpik(\tf{\nu}+\tf{x})\right) \notag\\
&\ge&(1-\varepsilon/2)\left(\sum_{i=1}^{k_1} C_i(\tf{v}^0)u^{-2}|x_i|^{\alpha_i+2(t_u^i)^{\beta_i}}
+\sum_{i=k_1+1}^{k} C_i(\tf{v}^0)u^{-2}|x_i|^{\alpha_i}\right)
\EQN
uniformly with respect to $\tf{\nu, \nu+x}\in[0,S]^k$.
Similarly, for sufficiently large $u$
\BQN
1-Cov\left(\Xjpik(\tf{\nu}),\Xjpik(\tf{\nu}+\tf{x})\right)%\notag\\
&\le&(1+\vn/2)\left(%\sum_{i=1}^{k_1} C_i(\tf{v}^0)u^{-2}(x_i)^{\alpha_i}
\sum_{i=1}^{k} C_i(\tf{v}^0)u^{-2}|x_i|^{\alpha_i}\right),
\EQN
uniformly with respect to $\tf{\nu, \nu+x}\in[0,S]^k$. The claim follows now \bb{by} the Slepian's inequality.\QED

\prooflem{lemsta2} The proofs of $(i)$ and $(ii)$ are similar, \Hr{therefore we present below only} the proof of $(i)$. Note that
\BQNY
\lim_{u\rw\IF} u^2(1-Cov(Y_{\vn,u}^{\tf{v}^0}(\tf{t}),Y_{\vn,u}^{\tf{v}^0}(\tf{s})))=(1-\vn)\sum_{i=1}^k \cJ{C_i(\tf{v}^0)}|t_i-s_i|^{\alpha_i}
\EQNY
uniformly with respect to $\tf{s,t}\in[0,S]^k$. Hence $(i)$ follows from \nelem{lemsta1}.\QED

\prooflem{lemdoubZ} Let
$$W_{\vn,u}(\tf{\nu,\nu'}):=\frac{1}{\sqrt{2}}\left(Z_{\vn,u}^{\tf{v}^0}(\tf{\nu})+\widetilde{Z}^{\tf{w}^0}_{\vn,u}(\tf{\nu'})\right),\ \ \ \tf{\nu,\nu'}\in[0,S]^{k}.$$
Since $\E{W_{\vn,u}(\tf{\nu,\nu'})}\equiv0$, $\E{(W_{\vn,u}(\tf{\nu,\nu'}))^2}\equiv1$, and
\BQNY
\lim_{u\rw\IF} u^2(1-Cov(W_{\vn,u}(\tf{\nu,\nu'}),W_{\vn,u}(\tf{\mu,\mu'})))=(1+\vn)\left(\sum_{i=1}^k \cJ{C_i(\tf{v}^0)}|\nu_i-\mu_i|^{\alpha_i}+\sum_{i=1}^k \cJ{C_i(\tf{w}^0)}|\nu'_i-\mu'_i|^{\alpha_i}\right)
\EQNY
uniformly with respect to $\tf{\nu,\mu,\nu',\mu'}\in[0,S]^k$,
it follows immediately from \nelem{lemsta1} that, as $u\rw\IF,$
\BQNY
&&\pk{\underset{\tf{\nu,\mu}\in[0,S]^{k}}\sup \frac{1}{\sqrt{2}}\left(Z_{\vn,u}^{\tf{v}^0}(\tf{\nu})+\widetilde{Z}^{\tf{w}^0}_{\vn,u}(\tf{\mu})\right)>u}\\
&&=\left(\prod_{i=1}^{k}\mH_{\alpha_i}\left[0,(C_i(\tf{v}^0) (1+\vn))^{1/\alpha_i}S\right]\right)\left(\prod_{i=1}^{k}\mH_{\alpha_i}\left[0,\cJ{C_i(\tf{w}^0)} (1+\vn))^{1/\alpha_i}S\right]\right)\Psi(u)(1+o(1))\\
&&\le \left(\prod_{i=1}^k \mH_{\alpha_i}[0,1](C_U^i(1+\vn))^{1/\alpha_i}\right)^2S^{2k}\Psi(u)(1+o(1)),
\EQNY
where in the last inequality we used the fact that $\mH_{\alpha_i}[0,R]\le \mH_{\alpha_i}[0,1]R,$ for any \bbb{$R>1$} (cf.\ Piterbarg (1996)), hence
\cH{the proof is complete.}  \QED

\prooflem{lemdoub} Since the proof of $(1)$ and $(2)$ are similar, we present next only the proof of $(1)$. Let
$$
Y_u(\tf{\nu,\nu'})=X_{1,u}(\tf{\nu})+X_{2,u}(\tf{\nu'}),
$$
where
$$X_{1,u}(\tf{\nu})=X\left(c_{p_1}^1+\frac{j_1 S+\nu_1}{u^{2/(\alpha_1(c_{p_1+1}^1))}},\cdots,c_{p_{k_1}}^{k_1}+\frac{j_1 S+\nu_{k_1}}{u^{2/(\alpha_{k_1}(c_{p_{k_1}+1}^{k_1}))}},
\tf{a}+\delta\kk+g_uS\ii +\triangle_0^{\tf{\nu}}\right)$$
and
$$X_{2,u}(\tf{\nu'})=X\left(c_{p'_1}^1+\frac{j'_1 S+\nu'_1}{u^{2/(\alpha_1(c_{p'_1+1}^1))}},\cdots,c_{p'_{k_1}}^{k_1}+\frac{j'_1 S+\nu'_{k_1}}{u^{2/(\alpha_{k_1}(c_{p'_{k_1}+1}^{k_1}))}},
\tf{a}+\delta\kk'+g_uS\iii +\triangle_0^{\tf{\nu'}}\right),$$
\cJ{ with $\triangle_0^{\tf{\nu'}}=g_u\prod_{i=k_1+1}^{k}[0,\nu'_i]$}. \cE{For any $u>0$, we have}
\BQNY
\pk{\underset{\tf{v}\in \Ajpik} %  A_{\vk{j,p},\ii,\kk}}
\sup X(\tf{v})>u, \underset{\tf{v'}\in \Ajpikp } %\prod_{i=1}^{k_1}B_{j'_i,p'_i}^i\times V_{\iii,\kk'}}
\sup X(\tf{v'})>u}\le\pk{\underset{\tf{\nu,\nu'}\in[0,S]^{k}}\sup Y_{u}(\tf{\nu,\nu'})>2u}.
\EQNY
We see from \eqref{eqcovvn} and \eqref{eqdis} that, for sufficiently large $u$, %uniformly for $(j_i,p_i),(j'_i,p'_i)\in \LL{1}^i, i=1,\cdots,k_1,(\ii,\kk),(\ii',\kk')\in\LL{2}$ satisfying $(\tf{j}, \tf{p})\preceq(\tf{j}', \tf{p}'),(\ii,\kk)\preceq(\ii',\kk')$,
\BQNY
Var(Y_u(\tf{\nu,\nu'}))=4-2(1-Cov(X_{1,u}(\tf{\nu}),X_{2,u}(\tf{\nu'})))>2.
\EQNY
It follows, for fixed $i=1,\cdots,k_1$, and $v_i\in B_{j_i,p_i}^i$, $v'_i\in B_{j'_i,p'_i}^i$, that
$|v_i-v'_i|\ge N_{j_i,p_i}^{j'_i,p'_i}\frac{S}{u^{2/(\alpha_i(c_{p_i+1}^i))}}$. Further, we have, for fixed $i=\gokkk$, $v_{k_1+i}\in \left[K_i\delta+\frac{I^{\kk}_i S}{u^{2/\alpha_{k_1+i}}},K_i\delta+\frac{(I^{\kk}_i+1) S}{u^{2/\alpha_{k_1+i}}}\right]$ and $v'_{k_1+i}\in \left[K'_i\delta+\frac{I'^{\kk'}_i S}{u^{2/\alpha_{k_1+i}}},K'_i\delta+\frac{(I'^{\kk'}_i+1) S}{u^{2/\alpha_{k_1+i}}}\right]$ that
$|v_{k_1+i}-v'_{k_1+i}|\ge \ N_{I^{\kk}_i,K_i}^{I'^{\kk'}_i,K'_i}\frac{S}{u^{2/\alpha_{k_1+i}}}$. Therefore, there exists some $\mathbb{C}_2>0$ such that for sufficiently large $u$
\BQNY
Var(Y_u(\tf{\nu,\nu'}))\le 4-\mathbb{C}_2\left(\sum_{i=1}^{k_1}\left(N_{j_i,p_i}^{j'_i,p'_i}\frac{S}{u^{2/\alpha_i(c_{p_i+1}^i)}}\right)^{\alpha_i(c_{p'_i+1}^i)}+\sum_{i=1}^{k_2}      \left(\ N_{I^{\kk}_i,K_i}^{I'^{\kk'}_i,K'_i}\frac{S}{u^{2/\alpha_{k_1+i}}}\right)^{\alpha_{k_1+i}}\right).
\EQNY
With the help of Lemma 4.4 of D\c{e}bicki and Kisowski (2008), we have, for some $\mathbb{C}_3>0$,
\BQNY
Var(Y_u(\tf{\nu,\nu'}))\le 4-\mathbb{C}_3\left(\sum_{i=1}^{k_1}\left(\sqrt{N_{j_i,p_i}^{j'_i,p'_i}}S\right)^{\alpha_i}+\sum_{i=1}^{k_2}
\left(N_{I^{\kk}_i,K_i}^{I'^{\kk'}_i,K'_i}S\right)^{\alpha_{k_1+i}}\right)u^{-2}\Hr{=:}H(S,u).
\EQNY
\cE{Consequently,} %We see from the last equation that
\BQNY
\pk{\underset{\tf{\nu,\nu'}\in[0,S]^{k}}\sup Y_{u}(\tf{\nu,\nu'})>2u}\le\pk{\underset{\tf{\nu,\nu'}\in[0,S]^{k}}\sup \overline{Y}_{u}(\tf{\nu,\nu'})>\frac{2u}{\sqrt{H(S,u)}}},
\EQNY
where  $\overline{Y}_{u}(\tf{\nu,\nu'})=Y_u(\tf{\nu,\nu'})/\sqrt{Var(Y_u(\tf{\nu,\nu'}))}$. Furthermore, following the argumentation analogous to that given in the proof of Lemma 6.3 in Piterbarg (1996) \cH{(see alternatively the proof of Lemma 4.5 in D\c{e}bicki and Kisowski (2008))},  for $\tf{\nu,\nu',\mu,\mu'}\in[0,S]^k$,
\BQNY
\E{(\overline{Y}_{u}(\tf{\nu,\nu'})-\overline{Y}_u(\tf{\mu,\mu'}))^2}&\le& 4\left(\E{(X_{1,u}(\tf{\nu})-X_{1,u}(\tf{\mu}))^2}+\E{(X_{2,u}(\tf{\nu'})-X_{2,u}(\tf{\mu'}))^2}\right)\\
&\le&\frac{1}{2}\left(\E{(Z^{\tf{v}^0}_{8,u}(\tf{\nu})-Z^{\tf{v}^0}_{8,u}(\tf{\mu}))^2}+\E{(\widetilde{Z}^{\tf{v'}^0}_{8,u}(\tf{\nu'})-\widetilde{Z}^{\tf{v'}^0}_{8,u}(\tf{\mu'}))^2}\right),
\EQNY
where \bbb{the GRF} \cE{$\widetilde{Z}^{\tf{v}'^0}_{8,u}$ is} independent of $Z_{8,u}^{\tf{v}^0}$, and has covariance structure \eqref{eqcovZ} with $\tf{v}^0$ replaced by $\tf{v'}^0$ \cJ{(chosen similarly as $\tf{v}^0$)}.
Next, \cH{by Slepian's} inequality \ee{(see e.g., Theorem C.1 of Piterbarg (1996))} and \nelem{lemdoubZ}, we obtain
\BQNY
&&\pk{\underset{\tf{\nu,\nu'}\in[0,S]^{k}}\sup \overline{Y}_{u}(\tf{\tf{\nu,\nu'}})>\frac{2u}{\sqrt{H(S,u)}}}\le
\pk{\underset{\tf{\nu,\nu'}\in[0,S]^{k}}\sup \frac{1}{\sqrt{2}}\left(Z^{\tf{v}^0}_{8,u}(\tf{\nu})+ \widetilde{Z}^{\tf{v'}^0}_{8,u}(\tf{\nu'})\right)>\frac{2u}{\sqrt{H(S,u)}}}\\
&&\qquad\le F_8S^{2k}\Psi\left(\frac{2u}{\sqrt{H(S,u)}}\right)\\
&&\qquad\le\mathbb{C}S^{2k}\expo{{-\mathbb{C}_1\left(\sum_{i=1}^{k_1}\left(\sqrt{N_{j_i,p_i}^{j'_i,p'_i}}S\right)^{\alpha_i}+\sum_{i=1}^{k_2}
\left(N_{I^{\kk}_i,K_i}^{I'^{\kk'}_i,K'_i}S\right)^{\alpha_{k_1+i}}\right)}}\Psi(u)%(1+o(1))
\EQNY
for $u$ sufficiently large.   %uniformly for $(j_i,p_i),(j'_i,p'_i)\in \LL{1}^i, i=1,\cdots,k_1,(\ii,\kk),(\ii',\kk')\in\LL{2}$.
Next, \cH{in order to} prove $(3)$ we apply the Borell theorem (e.g., Piterbarg (1996)). By \eqref{eqcov11} and \eqref{eqdis2}, we see that
\BQNY
\underset{\tf{v}\in \Ajpik,\tf{v'}\in \Ajpikp}\sup Var(X(\tf{v})+X(\tf{v'}))=4-2\underset{\tf{v}\in \Ajpik,\tf{v'}\in \Ajpikp}\inf(1-Cov(X(\tf{v}),X(\tf{v'})))< 4-2\lambda,
\EQNY
with some $\lambda\in(0,1)$. Further, there exists some $h>0$, such that
\BQNY
\pk{\underset{\tf{v}\in \Ajpik,\tf{v'}\in \Ajpikp}\sup X(\tf{v})+X(\tf{v'})>h}\le 2\pk{\underset{\tf{v}\in[0,T]^k}\sup X(\tf{v})>h/2}<\frac{1}{2}.
\EQNY
 Consequently, utilising Borell theorem, we obtain, for $u$ sufficiently large
 \BQNY
&&\pk{\underset{\tf{v}\in \Ajpik }\sup X(\tf{v})>u, \underset{\tf{v'}\in \Ajpikp }\sup X(\tf{v'})>u}\\
&&\le\pk{\underset{\tf{v}\in \Ajpik,\tf{v'}\in \Ajpikp}\sup X(\tf{v})+X(\tf{v'})>2u}\le 2\Psi\left(\frac{u-h/2}{\sqrt{1-\lambda/2}}\right)
\EQNY
establishing thus the claim.
\QED

\prooflem{lemtail1}
%In view of \nelem{lemtail} below, we only need to \cL{analyze} $\Pi(u)$.
Let $\vn\in(0,1)$ be \bbb{an } arbitrarily chosen
\bb{constant, and set $\OE:=1+ \ve$}. We first give the upper bound. Noting that $n_{i,p_i}=\lfloor \frac{c^i_{p_i+1}-c^i_{p_i}}{S}u^{2/\alpha_i(c_{p_i+1}^i)}\rfloor$, we derive that, as $u\rw\IF,$%,  for $u$ sufficiently large,
\BQNY
&&\Pi(u)\le\underset{(\ii,\kk)\in\UU{2}}{\underset{(j_i,p_i)\in \UU{1}^i, 1\le i\le k_1,}\sum}
\pk{\underset{\tf{v}\cE{\in \Ajpik}} %  A_{\vk{j,p},\ii,\kk}}
\sup X(\tf{v})>u}
\le\underset{(j_i,p_i)\in \UU{1}^i, 1\le i\le k_1}\sum\underset{\kk\in\mathcal{B}}\sum\underset{\ii\in\mathcal{A}_{\kk}}\sum\ \pk{\underset{\nu\in[0,S]^{k}}\sup Z_{\vn,u}^{\tf{v}^0}(\nu)>u}\\
&&\le\underset{p_i\le m_i, 1\le i\le k_1}\sum\ \underset{\kk\in\mathcal{B}}\sum\left(\prod_{i=1}^{k_1}\left(\frac{c_{p_i+1}^i-c_{p_i}^i}{S}u^{2/\left(\alpha_i(c_{p_i+1}^i)\right)}\right)N_{\kk}^+
\left(\prod_{i=1}^{k}\mH_{\alpha_i}[0,\cJ{C_i(\tf{v}^0)}\OE)^{1/\alpha_i}S]\right)\Psi(u)(1+o(1))\right)\\
&&=\underset{p_i\le m_i, 1\le i\le k_1}\sum\ \underset{\kk\in\mathcal{B}}\sum\Bigg(\frac{\prod_{i=1}^{k}\mH_{\alpha_i}[0,\cJ{C_i(\tf{v}^0)}\OE)^{1/\alpha_i}S]}{\prod_{i=1}^{k}(\cJ{C_i(\tf{v}^0)}\OE)^{1/\alpha_i}S)}
\left(\prod_{i=1}^{k}(\cJ{C_i(\tf{v}^0)}\OE)^{1/\alpha_i}S)\right)\frac{1}{S^{k_1}}\left(\prod_{i=1}^{k_1}\frac{u^{2/\alpha_i}}{(\ln u)^{1/\beta_i}}\right)\\
&&\times\frac{N_{\kk}^+\left(\prod_{i=k_1+1}^{k}(Su^{-2/\alpha_i})\right)}{\prod_{i=k_1+1}^{k}(Su^{-2/\alpha_i})}
\Psi(u)(1+o(1))
\prod_{i=1}^{k_1}\left((\ln u)^{1/\beta_i}(c_{p_i+1}^i-c_{p_i}^i)e^{\frac{2\left(\alpha_i-\alpha_i(c_{p_i+1}^i)\right)}{\alpha_i \alpha_i(c_{p_i+1}^i)}\ln u}\right)\Bigg)\\
&&\le\underset{p_i\le m_i, 1\le i\le k_1}\sum\ \underset{\kk\in\mathcal{B}}\sum\Bigg(\frac{\prod_{i=1}^{k}\mH_{\alpha_i}[0,\cJ{C_i(\tf{v}^0)}\OE)^{1/\alpha_i}S]}{\prod_{i=1}^{k}(\cJ{C_i(\tf{v}^0)}\OE)^{1/\alpha_i}S)}
\left(\prod_{i=1}^{k}(\cJ{C_i(\tf{v}^0)}\OE)^{1/\alpha_i})\right)
%\left(\frac{\prod_{i=1}^{k}u^{2/\alpha_i}}{\prod_{i=1}^{k_1}(\ln u)^{1/\beta_i}}\right)
\left(N_{\kk}^+\prod_{i=k_1+1}^{k}(Su^{-2/\alpha_i})\right)\\
&&\times\prod_{i=1}^{k_1}\left((\ln u)^{1/\beta_i}(c_{p_i+1}^i-c_{p_i}^i)
e^{-\frac{2(1-\vn)}{\alpha_i^2}\left((\ln u)^{1/\beta_i}c_{p_i+1}^i\right)^{\beta_i} }e^{\frac{2(1-\vn)}{\alpha_i^2}(\ln u)\left(c_{m_i+1}^i\right)^{\beta_i}\cJ{\abs{\ln(c_{m_i+1}^i)}^{-\delta_i} }  }
\right)\Bigg)\tLu\Psi(u)(1+o(1)),
\EQNY
\cL{where
$$
\tLu:=\frac{\prod_{i=1}^{k}u^{2/\alpha_i}}{\prod_{i=1}^{k_1}(\ln u)^{1/\beta_i}},$$
with $\prod_{i=m+1}^m(\cdot):=1, m\inn$.
}
 It follows that (see also D\c{e}bicki and Kisowski (2008))
\BQNY
\underset{S\rw\infty}\lim \frac{\prod_{i=1}^{k}\mH_{\alpha_i}[0,\cJ{C_i(\tf{v}^0)}\OE)^{1/\alpha_i}S]}{\prod_{i=1}^{k}(\cJ{C_i(\tf{v}^0)}\OE)^{1/\alpha_i}S)}=\prod_{i=1}^k\mH_{\alpha_i},
\ \ \
\lim_{u\rw\IF}e^{\frac{2(1-\vn)}{\alpha_i^2}(\ln u)(c_{m_i+1}^i)^{\beta_i}\cJ{\abs{\ln(c_{m_i+1}^i)}^{-\delta_i} }}=1,
\EQNY
%with $\mH_\alpha$ being the generalized Pickands' constant,
%\BQNY
%\lim_{u\rw\IF}\prod_{i=1}^{k_1}\cJ{C_i(\tf{v}^0)}(1+\vn))^{1/\alpha_i}=\prod_{i=1}^{k_1}(C_i(1+\vn))^{1/\alpha_i},\
%\EQNY
\BQNY
&&\underset{\delta\rw0}{\lim_{u\rw\IF}}\underset{\kk\in\mathcal{B}}\sum \prod_{i=\cJ{1}}^{k}(\cJ{C_i(\tf{v}^0)}\OE)^{1/\alpha_i}\left(N_{\kk}^+\prod_{i=k_1+1}^{k}(Su^{-2/\alpha_i})\right)
=\cJ{\int_{\tf{x}\in \prod_{i=1}^{k_1}\{t_i^0\}\times\prod_{i=k_1+1}^k[a_i,b_i]}\prod_{i=1}^{k}(C_i(\tf{x})\OE)^{1/\alpha_i}d{\tf{x}}}%\\
%&&\qquad\qquad= (1+\vn)^{\sum_{i=k_1+1}^k\frac{1}{\alpha_i} } \left(\prod_{i=k_1+1}^{k}\int_{a_i}^{b_i}(C_i(x))^{1/\alpha_i}dx\right),
\EQNY
and
\BQNY
&&\lim_{u\rw\IF}\underset{p_i\le m_i, 1\le i\le k_1}\sum\prod_{i=1}^{k_1}\left((\ln u)^{1/\beta_i}(c_{p_i+1}^i-c_{p_i}^i)
e^{-\frac{2(1-\vn)}{\alpha_i^2}\left((\ln u)^{1/\beta_i}c_{p_i+1}^i\right)^{\beta_i} } \right)=\int_{\R^{k_1}_+}e^{-\sum_{i=1}^{k_1}\frac{2(1-\vn)}{\alpha_i^2}x_i^{\beta_i}}d\tf{x}\\
&&\qquad\qquad=\prod_{i=1}^{k_1}\left(\frac{\alpha_i^2}{2(1-\vn)}\right)^{1/\beta_i} \frac{\Gamma(1/\beta_i)}{\beta_i}
\EQNY
since $\prod_{i=1}^{k_1}(\ln u)^{1/\beta_i}(c_{p_i+1}^i-c_{p_i}^i)\rw0$ and $(\ln u)^{1/\beta_i}c_{\bbb{m_i+1}}^i\rw\IF$, as $u\rw\IF$.
Consequently, the upper bound is given as
\BQNY
\Pi(u)&\le& %\left(\prod_{i=1}^{k_1}(C_i)^{1/\alpha_i}\right)
\OE^{\sum_{i=1}^k\frac{1}{\alpha_i} }\left(\prod_{i=1}^k\mH_{\alpha_i}\right) \left(\cJ{\int_{\tf{x}\in \prod_{i=1}^{k_1}\{t_i^0\}\times\prod_{i=k_1+1}^k[a_i,b_i]}\prod_{i=1}^{k}(C_i(\tf{x}))^{1/\alpha_i}d{\tf{x}}}\right)\\
&&\times\prod_{i=1}^{k_1}\left(\frac{\alpha_i^2}{2(1-\vn)}\right)^{1/\beta_i} \frac{\Gamma(1/\beta_i)}{\beta_i}
%\left(\frac{\prod_{i=1}^{k}u^{2/\alpha_i}}{\prod_{i=1}^{k_1}(\ln u)^{1/\beta_i}}\right)
\tLu\Psi(u)\cJ{(1+o(1))}
\EQNY
as $u\rw\IF.$ Next we derive the lower bound: using Bonferroni's inequality, we have
\BQNY
\Pi(u)&\ge&\underset{(\ii,\kk)\in\LL{2}}{\underset{(j_i,p_i)\in \LL{1}^i, 1\le i\le k_1,}\sum}\pk{\underset{\tf{v}\in \Ajpik} %  A_{\vk{j,p},\ii,\kk}}
\sup X(\tf{v})>u}\\
&&-2\underset{(\tf{j}, \tf{p})=(\tf{j}', \tf{p}')\ \mathrm{and}\ (\ii,\kk)\prec(\ii',\kk')}{\underset{(\tf{j}, \tf{p})\prec(\tf{j}', \tf{p}'),\mathrm{or}}{\underset{(j_i,p_i),(j'_i,p'_i)\in \LL{1}^i, 1\le i\le k_1,(\ii,\kk),(\ii',\kk')\in\LL{2}}\sum}}
\pk{\underset{\tf{v}\in \Ajpik} %  A_{\vk{j,p},\ii,\kk}}
\sup X(\tf{v})>u, \underset{\tf{v'}\in \Ajpikp} %\prod_{i=1}^{k_1}B_{j'_i,p'_i}^i\times V_{\ii',\kk'}}
\sup X(\tf{v'})>u}
%&\le&\underset{(j_i,p_i)\in \LL{1}^i, i=1,\cdots,k_1}\sum\underset{\ii}\sum\ \underset{\kk}\sum \pk{\underset{\nu\in[0,S]^{k}}\sup Z_{\vn,u}^{\tf{v}^0}(\nu)>u}\\
\EQNY
Similar arguments as \cH{in the derivation of} the upper bound yield, as $u\rw\IF,$%, for sufficiently large $u$,
\BQNY
&&\lim_{\delta\rw0,S\rw\IF}\underset{(\ii,\kk)\in\LL{2}}{\underset{(j_i,p_i)\in \LL{1}^i, 1\le i\le k_1,}\sum}
\pk{\underset{\tf{v}\in \Ajpik} %  A_{\vk{j,p},\ii,\kk}}
\sup X(\tf{v})>u}\\
&\ge&\lim_{\delta\rw0,S\rw\IF}\underset{(j_i,p_i)\in \LL{1}^i, 1\le i\le k_1}\sum\underset{(\ii,\kk)\in\LL{2}}\sum\ \pk{\underset{\nu\in[0,S]^{k}}\sup Y_{\vn,u}^{\tf{v}^0}(\nu)>u}\\
&\ge& %\left(\prod_{i=1}^{k_1}(C_i)^{1/\alpha_i}\right)
(1-\vn)^{\sum_{i=1}^k\frac{1}{\alpha_i} }\left(\prod_{i=1}^k\mH_{\alpha_i}\right) \left(\cJ{\int_{\tf{x}\in \prod_{i=1}^{k_1}\{t_i^0\}\times\prod_{i=k_1+1}^k[a_i,b_i]}\prod_{i=1}^{k}(C_i(\tf{x}))^{1/\alpha_i}d{\tf{x}}}\right)\\
&&\times\prod_{i=1}^{k_1}\left(\frac{\alpha_i^2}{2\OE}\right)^{1/\beta_i} \frac{\Gamma(1/\beta_i)}{\beta_i}
%\left(\frac{\prod_{i=1}^{k}u^{2/\alpha_i}}{\prod_{i=1}^{k_1}(\ln u)^{1/\beta_i}}\right)
\tLu\Psi(u)\cJ{(1+o(1))}.
\EQNY
Therefore, by letting  $\vn\rw0$, \cH{in order} to complete the proof, it is sufficient to show that
\BQN
&&\lim_{\delta\rw0,S\rw\IF}\lim_{u\rw\IF}\frac{\underset{(\tf{j}, \tf{p})=(\tf{j}', \tf{p}')\ \mathrm{and}\ (\ii,\kk)\prec(\ii',\kk')}{\underset{(\tf{j}, \tf{p})\prec(\tf{j}', \tf{p}'),\mathrm{or}}{\underset{(j_i,p_i),(j'_i,p'_i)\in \LL{1}^i, 1\le i\le k_1,(\ii,\kk),(\ii',\kk')\in\LL{2}}\sum}}
\pk{\underset{\tf{v}\in \Ajpik} %  A_{\vk{j,p},\ii,\kk}}
\sup X(\tf{v})>u, \underset{\tf{v'}\in \Ajpikp} %\prod_{i=1}^{k_1}B_{j'_i,p'_i}^i\times V_{\ii',\kk'}}
\sup X(\tf{v'})>u}}{\tLu\Psi(u)}\nonumber\\
%&&\qquad=\sum_{i=1}^3\lim_{u\rw\IF}\frac{\underset{(\tf{(j,p),(j',p')},(\ii,\kk),(\ii',\kk'))\in E_i}\sum
%\pk{\underset{\tf{v}\in \Ajpik}
%\sup X(\tf{v})>u, \underset{\tf{v'}\in \Ajpikp}
%\sup X(\tf{v'})>u}}{\tLu\Psi(u)}\nonumber\\
&&\qquad=\sum_{i=1}^3\lim_{\delta\rw0,S\rw\IF}\lim_{u\rw\IF}\frac{\Sigma^i_u}{\tLu\Psi(u)}=0, \label{eqdouble}
\EQN
where
$$
\Sigma^i_u:=\underset{(\tf{(j,p),(j',p')},(\ii,\kk),(\ii',\kk'))\in E_i}\sum
\pk{\underset{\tf{v}\in \Ajpik}
\sup X(\tf{v})>u, \underset{\tf{v'}\in \Ajpikp}
\sup X(\tf{v'})>u},\ i=1,2,3,
$$
with
\BQNY
&&E_\imath=\Big\{(\tf{(j,p),(j',p')},(\ii,\kk),(\ii',\kk')): \ \mq{conditions}\mq{of}\ (\imath)\ \mq{in}\ \mathrm{Lemma}\ 3.6\  \mq{are}\  \mathrm{satisfied}, \mq{and}\\
&&\qquad\qquad\qquad(\tf{j}, \tf{p})\prec(\tf{j}', \tf{p}'),\mathrm{or}\
 (\tf{j}, \tf{p})=(\tf{j}', \tf{p}')\ \mathrm{and}\ (\ii,\kk)\prec(\ii',\kk')\Big\}, \ \imath=1,2,3.
\EQNY
Eq. \eqref{eqdouble} follows from   \nelem{lemdoub}, and the  details are given in Appendix.\QED

\prooflem{lemtail} It is easy to see that the set $[0,T]^k/\prod_{i=1}^{k_1}[0,t^i_u] \times\prod_{i=k_1+1}^k[a_i,b_i]$ is the union of $2^{k_1}3^{k_2}-1$ sets of the form $\prod_{i=1}^{k_1}\Lambda_{i,u} \times\prod_{i=k_1+1}^k\Theta_i$, with
$$\Lambda_{i,u}=[0,t^i_u]\ \mathrm{or}\ [t^i_u, T], i=\gok, \quad \mathrm{and}\quad \Theta_i=[0,a_i] \ \mathrm{or}\ [a_i,b_i]\ \mathrm{or}\ [b_i,T],i=\gokk,$$
where at least one of $\{[t^i_u, T], i=\gok, [0,a_i], [b_i,T], i=\gokk\}$ appears. Since the other cases are similar, without loss of generality, it suffices to prove that
\BQNY
\pk{\underset{\tf{t}\in \prod_{i=1}^{k_1-1}[0,t^i_u]\times[t^{k_1}_u, T]\times\prod_{i=k_1+1}^{k-1}[a_i,b_i]\times[b_k,T]}{\sup}X(\tf{t})>u}=o(\Pi(u)).
\EQNY
We see that
\BQNY
&&\pk{\underset{\tf{t}\in \prod_{i=1}^{k_1-1}[0,t^i_u]\times[t^{k_1}_u, T]\times\prod_{i=k_1+1}^{k-1}[a_i,b_i]\times[b_k,T]}{\sup}X(\tf{t})>u}\\
&&\le \pk{\underset{\tf{t}\in \prod_{i=1}^{k_1-1}[0,t^i_u]\times[t^{k_1}_u, T]\times\prod_{i=k_1+1}^{k-1}[a_i,b_i]\times[b_k,b_k+t^k_u]}{\sup}X(\tf{t})>u}\\
&&+\pk{\underset{\tf{t}\in \prod_{i=1}^{k_1-1}[0,t^i_u]\times[t^{k_1}_u, T]\times\prod_{i=k_1+1}^{k-1}[a_i,b_i]\times[b_k+t^k_u,T]}{\sup}X(\tf{t})>u}
\EQNY
It is sufficient to analyze the first probability on the right-hand side of the last inequality since the analysis of the  second one is similar.
It is derived that
\BQN
\theta(u)&:=&\pk{\underset{\tf{t}\in \prod_{i=1}^{k_1-1}[0,t^i_u]\times[t^{k_1}_u, T]\times\prod_{i=k_1+1}^{k-1}[a_i,b_i]\times[b_k,b_k+t^k_u]}{\sup}X(\tf{t})>u}\nonumber\\
&\le&\underset{(j_i,p_i)\in \UU{1}^i, i=1,\cdots,k_1-1,k,(\ii,\kk)\in\UU{2}'}\sum\pk{\underset{\tf{v}\in\prod_{i=1}^{k_1-1}B_{j_i,p_i}^i\times[t_u^{k_1},T]\times W_{\ii,\kk}\times\left( b_k+B^{k}_{j_k,p_k}\right)}\sup X(\tf{v})>u},\label{eqtail1}
%&\le&\underset{(j_i,p_i)\in \UU{1}^i, i=1,\cdots,k_1}\sum\underset{\kk}\sum\underset{\ii}\sum\ \pk{\underset{\nu\in[0,S]^{k}}\sup Z_{\vn,u}^{\tf{v}^0}(\nu)>u}\\
\EQN
where $\kk=(K_1,\cdots,K_{k_2-1})\in\mathbb{Z}^{k_2-1},$ $\ii=(I^{\kk}_1,\cdots, I^{\kk}_{k_2-1})\in\mathbb{Z}^{k_2-1}$, and $B^{k}_{j_k,p_k},$ $\mathcal{U}'_2$ and $W_{\ii,\kk}$ are defined similarly as $B^{k_1}_{j_{k_1},p_{k_1}}$, $\mathcal{U}_2$ and $V_{\ii,\kk}$, respectively.

For any fixed $j_i,p_i,i=1,\cdots, k_1,k, \ii, \kk$ \cJ{such that $(j_i,p_i)\in \UU{1}^i, i=1,2,\cdots,k_1-1,k$ and  $(\ii,\kk)\in\UU{2}'$},  consider the  GRF $X(\tf{v}):=X(v_1,\cdots,v_{k})$ \cE{on \cH{the} set}
$$\mathcal{A}_{\tf{jp,\ii,\kk}}:=\prod_{i=1}^{k_1-1}B_{j_i,p_i}^i\times[t_u^{k_1},T]\times W_{\ii,\kk}\times\left( b_k+B^{k}_{j_k,p_k}\right).$$
\cE{For notational simplicity \Hr{write next $X_{k_1,u}(\tf{\nu})$ instead of}}
$$X\left(c_{p_1}^1+\frac{j_1 S+\nu_1}{u^{2/(\alpha_1(c_{p_1+1}^1))}},\cdots,c_{p_{k_1-1}}^{k_1-1}+\frac{j_{k_1-1} S+\nu_{k_1-1}}{u^{2/(\alpha_{k_1-1}(c_{p_{k_1-1}+1}^{k_1-1}))}},\nu_{k_1},
\tf{a'}+\delta\kk+g'_uS\ii +{\triangle'}_0^{\tf{\nu}},b_k+c_{p_k}^k+\frac{j_k S+\nu_k}{u^{2/(\alpha_k(c_{p_k+1}^k))}}\right),$$
 where ${\triangle'}_0^{\tf{\nu}}=g'_u\prod_{i=k_1+1}^{k-1}[0,\nu_i]$, $\tf{a'}=(a_{k_1+1},\cdots,a_{k-1})$ and $g'_u$ is defined in a similar way as $g_u$ (see \eqref{eqgu}). It follows that
\BQN
\underset{\tf{v}\in \mathcal{A}_{\tf{jp,\ii,\kk}}} %\prod_{i=1}^{k_1-1}B_{j_i,p_i}^i\times[t_u^{k_1},T]\times W_{\ii,\kk}\times\left( b_k+B^{k}_{j_k,p_k}\right)}
\sup X(\tf{v})\ \overset{d}=\  \underset{\tf{\nu}\in[0,S]^{k_1-1}\times[t_u^{k_1},T]\times[0,S]^{k_2}}\sup X_{k_1,u}(\tf{\nu}).
\EQN
 Let $b_{k_1,u}=u^{-2/{\left(\alpha_{k_1}+\frac{3}{4}(t^{k_1}_u)^{\beta_{k_1}}\right)}}$, and fix
 $\tf{v}^0\in\prod_{i=1}^{k_1-1}A_{p_i}^i\times[0,T]\times \delta_{\kk}\times(b_k+A_{p_k}^k)$ with $A_{p_i}^i,\delta_{\kk}$ defined similarly as before (the only difference is the dimension).
In view of the proof of \eqref{slep2}, there exists a constant $\mathbb{C}_0$ such that, for sufficiently large $u$
\BQNY
1-Cov(X_{k_1,u}(\tf{\nu}),X_{k_1,u}(\tf{\nu}+\tf{x}))\le 1- e^{-\frac{3}{2}\sum_{i=1,i\neq k_1}^{k} C_i(\tf{v}^0)u^{-2}|x_i|^{\alpha_i}
-\mathbb{C}_0|x_{k_1}|^{\alpha_{k_1}+\frac{3}{4}(t_u^{k_1})^{\beta_{k_1}}}}
\EQNY
uniformly with respect to $\tf{\nu, \nu}+\tf{x}\in[0,S]^{k_1-1}\times[t_u^{k_1},T]\times[0,S]^{k_2}$ such that $|x_{k_1}|\le b_{k_1,u}$. Let $\{\widetilde{Z}_{u}^{\tf{v}^0}(\tf{t}), \tf{t}\in[0,S]^{k_1-1}\times[t_u^{k_1},T]\times[0,S]^{k_2}\}$, $u>0$, %with $x_{k_1}\le b_{k_1,u}$,
  \cH{be} a family of centered stationary  GRF's   such that
 \BQNY
Cov(\widetilde{Z}_{u}^{\tf{v}^0}(\tf{\nu}),\widetilde{Z}_{u}^{\tf{v}^0}(\tf{\nu}+\tf{x}))=e^{-\frac{3}{2}\sum_{i=1,i\neq k_1}^{k} C_i(\tf{v}^0)u^{-2}|x_i|^{\alpha_i}
-\mathbb{C}_0|x_{k_1}|^{\alpha_{k_1}+\frac{3}{4}(t_u^{k_1})^{\beta_{k_1}}}}
\EQNY
for $u$ such that $\alpha_{k_1}+\frac{3}{4}(t_u^{k_1})^{\beta_{k_1}}\le 2$, and $\tf{\nu, \nu}+\tf{x}\in[0,S]^{k_1-1}\times[t_u^{k_1},T]\times[0,S]^{k_2}$.
In view of the Slepian's inequality, continuing \eqref{eqtail1} we get, as $u\rw\IF$
\BQNY
\theta(u)%\pk{\underset{\tf{t}\in \prod_{i=1}^{k_1-1}[0,t^i_u]\times[t^{k_1}_u, T]\times\prod_{i=k_1+1}^{k-1}[a_i,b_i]\times[b_k,b_k+t^k_u]}{\sup}X(\tf{t})>u}\nonumber\\
&\le&\underset{(j_i,p_i)\in \UU{1}^i, i=1,\cdots,k_1-1,k,(\ii,\kk)\in\UU{2}'}\sum\pk{\underset{\tf{v}\in \mathcal{A}_{\tf{jp,\ii,\kk}}} %\prod_{i=1}^{k_1-1}B_{j_i,p_i}^i\times[t_u^{k_1},T]\times W_{\ii,\kk}\times\left( b_k+B^{k}_{j_k,p_k}\right)}
\sup X(\tf{v})>u}\\
&\le&\underset{(j_i,p_i)\in \UU{1}^i, i=1,\cdots,k_1-1,k}\sum\underset{(\ii,\kk)\in\UU{2}'}\sum\ \sum_{l=0}^{\lfloor T(b_{k_1,u})^{-1}\rfloor+1}\ \pk{\underset{\nu\in[0,S]^{k_1-1}\times[lb_{k_1,u},(l+1)b_{k_1,u}]\times[0,S]^{k_2}}\sup X_{k_1,u}(\tf{\nu})>u}\\
&\le&(\lfloor T(b_{k_1,u})^{-1}\rfloor+2)\underset{(j_i,p_i)\in \UU{1}^i, i=1,\cdots,k_1-1,k}\sum\underset{(\ii,\kk)\in\UU{2}'}\sum\ \pk{\underset{\tf{\nu}\in[0,S]^{k_1-1}\times[0,b_{k_1,u}]\times[0,S]^{k_2}}\sup \widetilde{Z}_{u}^{\tf{v}^0}(\tf{\nu})>u}\\
&\le&\left(u^{2/\alpha_{k_1}}(\ln u)^{-\frac{4}{3\beta_{k_1}}}T+2\right)\underset{(j_i,p_i)\in \UU{1}^i, i=1,\cdots,k_1-1,k}\sum\underset{(\ii,\kk)\in\UU{2}'}\sum\ \pk{\underset{\tf{\nu}\in[0,S]^{k_1-1}\times[0,b_{k_1,u}]\times[0,S]^{k_2}}\sup \widetilde{Z}_{u}^{\tf{v}^0}(\tf{\nu})>u},
\EQNY
where in the last inequality we used that $(b_{k_1,u})^{-1}\le u^{2/\alpha_{k_1}}(\ln u)^{-\frac{4}{3\beta_{k_1}}}$ given in D\c{e}bicki and Kisowski (2008).
Furthermore, it follows from \nelem{lemsta1} that, as $u\rw\IF$,
\BQNY
&&\pk{\underset{\nu\in[0,S]^{k_1-1}\times[0,b_{k_1,u}]\times[0,S]^{k_2}}\sup \widetilde{Z}_{u}^{\tf{v}^0}(\tf{\nu})>u}\\
&&=\left(\prod_{i=1}^{k_1-1}\mH_{\alpha_i}\left[0,\left(\frac{3}{2}C_i(\tf{v}^0)\right)^{1/\alpha_i}S\right]\times\mH_{\alpha_{k_1}}[0,\mathbb{C}_0^{1/\alpha_{k_1}}]
\times\prod_{i=k_1+1}^{k}\mH_{\alpha_i}\left[0,\left(\frac{3}{2}C_i(\tf{v}^0)\right)^{1/\alpha_i}S\right]\right)\Psi(u)(1+o(1))\\
&&\le \mathbb{C}_3\prod_{i=1}^k\mH_{\alpha_i}[0,1]S^{k-1}\Psi(u)(1+o(1))
\EQNY
\COM{
\cE{Next, Lemma 6.2 of Piterbarg (1996) implies}
\BQNY
\mH_{(L,\alpha)}\left[\prod_{i=1}^{k_1-1}\left[0,\left(\frac{3}{2}C_i(v_i^0)\right)^{1/\alpha_i}S\right]\times[0,1]
\times\prod_{i=k_1+1}^{k}\left[0,\left(\frac{3}{2}C_i(v_i^0)\right)^{1/\alpha_i}S\right]\right]\le
\mathbb{C}_3\prod_{i=1,i\neq k_1}^k \left(\left(\frac{3}{2}C_i(v_i^0)\right)^{1/\alpha_i}S\right)
\EQNY
}
for some positive constant $\mathbb{C}_3$. Consequently,  similar arguments as in the proof of the upper bound in \netheo{mainthm} \bb{implies}
\BQNY
\theta(u)
&\le& \mathbb{C}_4T%T\left(\prod_{i=1}^{k_1-1}(C_i)^{1/\alpha_i}\right)(C_{k}(b_k))^{1/\alpha_k}\left(\frac{3}{2}\right)^{\sum_{i=1,i\neqk_1}^{k}\frac{1}{\alpha_i} }
\left(\prod_{i=k_1+1}^{k-1}(b_i- a_i)\right)\left(\prod_{i=1}^{k_1-1} \frac{(\alpha_i)^{2/\beta_i} \Gamma(1/\beta_i)}{\beta_i} \frac{(\alpha_k)^{2/\beta_k} \Gamma(1/\beta_k)}{\beta_k} \right)%\int_{\R^{k_1}_+}e^{-\sum_{i=1}^{k_1-1}\frac{x_i^{\beta_i}}{\alpha_i^2}-\frac{x_k^{\beta_k}}{\alpha_k^2}}d\tf{x}
\left(\frac{\prod_{i=1}^{k}u^{2/\alpha_i}}{\prod_{i=1}^{k_1-1}(\ln u)^{1/\beta_i}}\right)(\ln u)^{-\frac{4}{3\beta_{k_1}}-\frac{1}{\beta_k}}\Psi(u)\\
&=&o(\Pi(u))
\EQNY
as $u\rw\IF$, and thus  the proof is complete.
 \QED

\section{Appendix}

\prooflem{lemsta1} Using the \cH{classical approach} (see e.g., Piterbarg (1996)) we have for $\cL{u>0}$ %(set \bbb{$\vk{0}=(0 \ldot 0)\inr^k$})
\BQN
\pk{\underset{\tf{t}\in \tf{D}}\sup X_u(\tf{t})>u}
&=&\frac{1}{\sqrt{2\pi}u}e^{-\frac{u^2}{2}}\int_{-\IF}^\IF e^{z-\frac{z^2}{2u^2}}\pk{\underset{\tf{t}\in\tf{D}}\sup X_u(\tf{t})>u|X_u(\tf{0})=u-\frac{z}{u}}\cH{dz}.\label{eqXu}
\EQN
It follows that, for any $u> 0$
\BQNY
\left\{X_u(\tf{t})|X_u(\tf{0})=u-\frac{z}{u}, \tf{t}\in\tf{D}\right\}\ \mathrm{and}\ \left\{X_u(\tf{t})-r_u(\tf{t},\tf{0})X_u(\tf{0})+r_u(\tf{t},\tf{0})\left(u-\frac{z}{u}\right), \tf{t}\in\tf{D}\right\}
\EQNY
have the same distribution (cf.\ Aldler and Taylor (2007) %or Tabi\'{s} 2011),
from which we see that
\BQNY
\pk{\underset{\tf{t}\in\tf{D}}\sup X_u(\tf{t})>u|X_u(\tf{0})=u-\frac{z}{u}}
=\pk{\underset{\tf{t}\in\tf{D}}\sup \left(\zeta_u(\tf{t})-u^2 (1-r_u(\tf{t},\tf{0}))+z(1-r_u(\tf{t},\tf{0}))\right)>z},
\EQNY
with $\left\{\zeta_u(\tf{t})=u(X_u(\tf{t})-r_u(\tf{t},\tf{0})X_u(\tf{0})), \tf{t}\in\tf{D}\right\}$.
By \eqref{eqlimru}%, we obtain that
\BQNY
\lim_{u\rw\IF}(u^2 (1-r_u(\tf{t},\tf{0}))-z(1-r_u(\tf{t},\tf{0})))=|\tf{t}|_{\tf{\alpha}}
\EQNY
uniformly with respect to $\tf{t}\in\tf{D}$ for any $z\in\R$.

Next we show that $\zeta_u, u >0$ \bb{converges} weakly to $\widetilde{B}_{\tf{\alpha}}$ in $\mathcal{C}(\tf{D})$ as $u\rightarrow\IF$. To this end, we need to show (e.g., Wichura (1969) or Neuhaus (1971)):

$i)$ finite-dimensional distributions of $\zeta_u$ converge in distribution to those of $\widetilde{B}_{\tf{\alpha}}$ as $u\rightarrow\IF$

$ii)$  tightness, i.e., %of $\{\zeta_u, u\ge 0\}$ in $\mathcal{C}(\tf{D})$.
for any $\eta>0$
\BQNY
\lim_{\delta\rightarrow0}\limsup_{u\rightarrow\IF}\pk{\underset{\max_{1\le i\le k}|s_i-t_i|<\delta}{\underset{\tf{s},\tf{t}\in\tf{D}}\sup}|\zeta_u(\tf{t})-\zeta_u(\tf{s})|>\eta}=0.
\EQNY
\bb{First note that the} increments of the centered  GRF   $\{\zeta_u(\tf{t}), \tf{t}\in\tf{D}\}$ have the following property
 \BQN
 \lim_{u\rightarrow\IF}Var(\zeta_u(\tf{t})-\zeta_u(\tf{s}))&=& \lim_{u\rightarrow\IF}\E{(\zeta_u(\tf{t})-\zeta_u(\tf{s}))^2}\nonumber\\
& = &\lim_{u\rightarrow\IF}(2u^2(1-r_u(\tf{t},\tf{s}))-u^2(r_u(\tf{t},\tf{0})-r_u(\tf{s},\tf{0}))^2)\nonumber\\
&=&2{|\tf{t}-\tf{s}|}_{\alpha}\nonumber\\
&=&Var(\widetilde{B}_{\tf{\alpha}}(\tf{t})-\widetilde{B}_{\tf{\alpha}}(\tf{s})).\label{eqvarZu}
 \EQN
\bb{Furthermore, the above holds} uniformly with respect to  $\tf{t}, \tf{s}\in \tf{D}$, implying $i)$. In order to prove the tightness, we use \cH{a} similar approach as in Dieker (2005) and D\c{e}bicki et al. (2012). We start by defining, for  fixed $u>0$, a metric $d_u$ on $\R_+^k$ as
\BQNY
d_u(\tf{s},\tf{t})=\sqrt{\E{(\zeta_u(\tf{t})-\zeta_u(\tf{s}))^2}}.
\EQNY
Further %, for any $d_u$-compact set $\tf{D}'\in \R_+^d$,
 write
\BQNY
B_{d_u}(t,u,\vartheta):=\{\tf{s}\in\R_+^k: d_u(\tf{s},\tf{t})\le\vartheta\}
\EQNY
for the $d_u$-ball centered at $\tf{t}\in\R_+^k$ and of radius $\vartheta$,
and let
\BQNY
H_{d_u}(\tf{D}',u, \vartheta):=\ln (N'(\tf{D}',u, \vartheta)),
\EQNY
with $N'(\tf{D}',u, \vartheta)$ being the smallest number of such balls that cover $\tf{D}'$, a compact set in $\R_+^k$. Here $H_{d_u}(\tf{D}',u, \vartheta)$ is called {\it (metric) entropy} for $\tf{D}'$ induced by $d_u$. See
%\cH{Lifshits (1995, 2012)} or
 Adler and Taylor (2007) for more detail on metric entropy.

We see from \eqref{eqvarZu} that, for $u$ sufficiently large, there exists some constant \bbb{$\bbb{\mathbb{C}_0}$} such that
\BQN
d_u(\tf{s},\tf{t})\le \bbb{\mathbb{C}_0}\sqrt{|\tf{s}-\tf{t}|_\alpha}\le k \bbb{\mathbb{C}_0} \delta^{\frac{\tf{\alpha}}{2}},\label{eqdu}
\EQN
if $\max_{1\le i\le k}|s_i-t_i|<\delta<1$, where
$\alpha:=\min_{1\le i\le k} \alpha_i$.
By utilising Corollary 1.3.4 of Adler and Taylor (2007), it follows that there exists some universal constant $Q_0>0$ such that, for any $0<\delta<1$,
\BQNY
\pk{\underset{\max_{1\le i\le k}|s_i-t_i|<\delta}{\underset{\tf{s},\tf{t}\in\tf{D}}\sup}|\zeta_u(\tf{t})-\zeta_u(\tf{s})|>\eta}&\le& \pk{\underset{d_u(\tf{s},\tf{t})<k \bbb{\mathbb{C}_0}\delta^{\frac{\alpha}{2}}}{\underset{\tf{s},\tf{t}\in\tf{D}}\sup}|\zeta_u(\tf{t})-\zeta_u(\tf{s})|>\eta}\\
&\le&\frac{Q_0}{\eta}\int_{0}^{k\bbb{\mathbb{C}_0} \delta^{\frac{\alpha}{2}}}\sqrt{H_{d_u}([0,R]^k,u,\vartheta)}d\vartheta,
\EQNY
with $R<\IF$ being a sufficiently large constant.
Define, for $\tf{t},\tf{s}\in\R_+^k$, a semimetric
\BQNY
\tilde{d}(\tf{t},\tf{s})=\bbb{\mathbb{C}_0} \sqrt{|\tf{s}-\tf{t}|_\alpha}.
\EQNY
Thanks to \eqref{eqdu}  it follows that, for sufficiently large $u$ and small $\vartheta$,
\BQNY
H_{d_u}([0,R]^k,u,\vartheta)\le H_{\tilde{d}}([0,R]^k,u,\vartheta)\le k\ln\left(\frac{R}{\left(\frac{\vartheta^2}{k \bbb{\mathbb{C}_0}^2}\right)^{\frac{1}{\alpha}}}+1\right)\le \bbb{\bbb{\mathbb{C}_1}}\ln\left(\frac{1}{\vartheta}\right),
\EQNY
for some positive constant $\bbb{\mathbb{C}_1}$, with $H_{\tilde{d}}([0,R]^k,u,\vartheta)$ being the entropy induced by $\tilde{d}$.

 Consequently, we have that
\BQNY
\lim_{\delta\rightarrow0}\limsup_{u\rightarrow\IF}\pk{\underset{\max_{1\le i\le k}|s_i-t_i|<\delta}{\underset{\tf{s},\tf{t}\in\tf{D}}\sup}|\zeta_u(\tf{t})-\zeta_u(\tf{s})|>\eta}\le \lim_{\delta\rw0}\frac{Q_0\sqrt{\bbb{\mathbb{C}_1}}}{\eta}\int_{\frac{1}{k\bbb{\mathbb{C}_0}}\delta^{-\frac{\alpha}{2}}}^\IF\frac{\sqrt{\ln \vartheta}}{\vartheta^2}d\vartheta=0,
\EQNY
establishing the claim $ii)$. Moreover, since the functional $\sup_{\tf{t}\in\tf{D}}f(\tf{t})$ is continuous on $\mathcal{C}(\tf{D})$, we conclude, for any $z\in\R$, that
\BQNY
\lim_{u\rw\IF}\pk{\underset{\tf{t}\in\tf{D}}\sup X_u(\tf{t})>u|X_u(\tf{0})=u-\frac{z}{u}}
=\pk{\underset{\tf{t}\in\tf{D}}\sup (\widetilde{B}_{\tf{\alpha}}(\tf{t})-|\tf{t}|_{\tf{\alpha}})>z}.
\EQNY
In order to use dominate convergence theorem to the integral in \eqref{eqXu} when taking limit in $u$, we need a uniform (in $u$ large enough) upper bound of
\BQNY
P_u(z):=\pk{\underset{\tf{t}\in\tf{D}}\sup \left(\zeta_u(\tf{t})-u^2 (1-r_u(\tf{t},\tf{0}))+z(1-r_u(\tf{t},\tf{0}))\right)>z}
\EQNY
for $z>0$ sufficiently large. It follows that, for $u$ sufficiently large,
\BQN
P_u(z)&\le& \pk{\underset{\tf{t}\in\tf{D}}\sup \zeta_u(\tf{t})+\underset{\tf{t}\in\tf{D}}\sup (1-r_u(\tf{t},\tf{0}))z>z}\nonumber\\
&\le& \pk{\underset{\tf{t}\in\tf{D}}\sup \zeta_u(\tf{t})>(1-\varrho_0)z}
\label{equperpu}
\EQN
for some $\varrho_0\in(0,1)$.
Further, we see from \eqref{eqvarZu} that, for sufficiently large $u$, there exists some positive constant \bb{$\bbb{\mathbb{C}_2}$} such that
\BQNY
Var(\zeta_u(\tf{t})-\zeta_u(\tf{s}))\le \bbb{\mathbb{C}_2} Var(\widetilde{B}_{\tf{\alpha}}(\tf{t})-\widetilde{B}_{\tf{\alpha}}(\tf{s}))
\EQNY
for all $\tf{s},\tf{t}\in \tf{D}$, implying, by Sudakov-Fernique inequality (e.g., Adler and Taylor (2007))
\BQN
\E{\sup_{\tf{t}\in\tf{D}}\zeta_u(\tf{t})}\le \sqrt{\bbb{\mathbb{C}_2}}\E{\sup_{\tf{t}\in\tf{D}}\widetilde{B}_{\tf{\alpha}}(\tf{t})}:=\cJ{U_0}.\label{eqborell1}
\EQN
The constant $\cJ{U_0}$ is finite, which follows thanks to Theorem 2.1.1 of Adler and Taylor (2007). %Therefore, we derive, for sufficiently large $u$, that
%\BQN
%\E{\underset{\tf{t}\in\tf{D}}\sup \left(\zeta_u(\tf{t})-u^2 (1-r_u(\tf{t},\tf{0}))+w(1-r_u(\tf{t},\tf{0}))\right)}\le \cJ{U_0}+\varrho_0 w
%\EQN
%for some small positive constant $\varrho_0$.
Moreover, %it is easy to see that
\BQN
\sup_{\tf{t}\in\tf{D}}Var(\zeta_u(\tf{t}))\le\sigma_{\tf{D}}^2:=\bbb{\mathbb{C}_2}\sup_{\tf{t}\in\tf{D}}Var(\widetilde{B}_{\tf{\alpha}}(\tf{t}))=2\bbb{\mathbb{C}_2}\sup_{\tf{t}\in\tf{D}} |\tf{t}|_{\tf{\alpha}}<\IF. \label{eqborell2}
\EQN
With the help of \eqref{equperpu}, \eqref{eqborell1} and \eqref{eqborell2}, Borell-TIS inequality (Theorem 2.1.1 of Adler and Taylor (2007)) gives, for any $z>\frac{\cJ{U_0}}{1-\varrho_0}$ and $u$ sufficiently large,
\BQNY
P_u(z)\le\pk{\underset{\tf{t}\in\tf{D}}\sup \zeta_u(\tf{t})>(1-\varrho_0)z}\le \exp\left(-\frac{((1-\varrho_0)z-\cJ{U_0})^2}{2\sigma_{\tf{D}}^2}\right).
\EQNY
Applying dominate convergence theorem to the integral in \eqref{eqXu}, we conclude that
\BQNY
\lim_{u\rw\IF}\int_{-\IF}^\IF e^{z-\frac{z^2}{2u^2}}P_u(z)dz=\E{\exp\left(\underset{\tf{t}\in\tf{D}}\sup (\widetilde{B}_{\tf{\alpha}}(\tf{t})-|\tf{t}|_{\tf{\alpha}})\right)},
\EQNY
\cH{thus the proof is completed}. \QED

{\textsc{Proof of}} Eq.\ \eqref{eqdouble} According to \nelem{lemdoub}, the three parts of the double-sum in \eqref{eqdouble} can be  estimated in different ways. It follows from \eqref{eqdoub3} that%$(3)$ of \nelem{lemdoub} that
\BQNY
\lim_{u\rw\IF}\frac{\Sigma^3_u}{\tLu\Psi(u)}\le \lim_{u\rw\IF}\frac{2\Psi\left(\frac{u-h/2}{\sqrt{1-\lambda/2}}\right) \left(\underset{(\ii,\kk)\in\LL{2},(j_i,p_i)\in \LL{1}^i, 1\le i\le k_1}\sum 1\right)^2}{\tLu\Psi(u)}=0,
\EQNY
where the sum in the middle term can be estimated using the same arguments as the upper bound in \netheo{mainthm}. Next, for sake of simplicity, we only give the estimates of the first two sums for $k_1=k_2=1$, since the general cases ($k_1,k_2$ are arbitrary integers) follow from   similar arguments. For the first sum, we derive, using  \eqref{eqdoub1} that, for $u$ sufficiently large %$u\rw\IF,$
\BQNY
&&\Sigma_u^1\le\underset{(I^{\kk}_1,K_1)\in\LL{2},(j_1,p_1)\in \LL{1}^1}\sum \Bigg(4\!\!\!\underset{(j_1,p_1)\prec(j'_1,p'_1) \mq{and} N_{j_1,p_1}^{j'_1,p'_1}>0}{\underset{(j'_1,p'_1)\in \LL{1}^1}\sum} \  %(j_1,p_1)\in \LL{1}^1,
\underset{N_{I^{\kk}_1,K_1}^{I'^{\kk'}_1,K'_1}\ge0}{\underset{(I'^{\kk'}_1,K'_1)\in\LL{2}}\sum}
\mathbb{C}S^{4}\exp\Big(-\mathbb{C}_1\Big((N_{j_1,p_1}^{j'_1,p'_1})^{\alpha_1/2}S^{\alpha_1}\\
&&+\left(N_{I^{\kk}_1,K_1}^{I'^{\kk'}_1,K'_1}\right)^{\alpha_2}S^{\alpha_{2}}\Big)\Big)
+ 2\underset{N_{I^{\kk}_1,K_1}^{I'^{\kk'}_1,K'_1}>0}{\underset{(I'^{\kk'}_1,K'_1)\in\LL{2}}\sum} \mathbb{C}S^{4}\expo{-\mathbb{C}_1
\left(N_{I^{\kk}_1,K_1}^{I'^{\kk'}_1,K'_1}\right)^{\alpha_2}S^{\alpha_{2}}}\Bigg)\Psi(u)\\
&&\le 4\mathbb{C}S^{4}\underset{(I^{\kk}_1,K_1)\in\LL{2},(j_1,p_1)\in \LL{1}^1}\sum \left(\left(\sum_{n_1\ge1}e^{-\mathbb{C}_1(\sqrt{n_1} S)^{\alpha_1}}\right)
\left(\sum_{n_2\ge0}e^{-\mathbb{C}_1({n_2} S)^{\alpha_2}}\right)+ \left(\sum_{n_3\ge1}e^{-\mathbb{C}_1({n_3} S)^{\alpha_2}}\right) \right)\Psi(u)\\
&&\le \mathbb{C}'S^{4}\underset{(I^{\kk}_1,K_1)\in\LL{2},(j_1,p_1)\in \LL{1}^1}\sum \left(e^{-\mathbb{C}'_1 S^{\alpha_1}}
\left(1+e^{-\mathbb{C}''_2  S^{\alpha_1}}\right)+ e^{-\mathbb{C}'''_2 S^{\alpha_1}} \right)\Psi(u),
\EQNY
for suitably chosen constants. % $\mathbb{C}',\mathbb{C}'_1, \mathrm{C}''_1, \mathrm{C}'''_1$.
This, combined with the estimate of the last sum in the above formula, yields that
\BQN
\lim_{S\rw\IF}\lim_{u\rw\IF}\frac{\Sigma^1_u}{\tLu\Psi(u)}=0.\label{eqSig2}
\EQN
Lastly we estimate the second sum.
 \bb{According} to $(2)$ of \nelem{lemdoub}, the sum $\Sigma_u^2$ can be divided into four parts, denoted by  $\Sigma_{i_1^1,u}^2, \Sigma_{i_2^1,u}^2, \Sigma_{i_1^2,u}^2$ and $ \Sigma_{i_2^2,u}^2$, respectively.
\cH{ Applying} \eqref{eqdoub2}, \nelem{lemslep1} and \nelem{lemsta2} we find that, for $u$ large enough, %as $u\rw\IF$
\BQNY
&&\Sigma_{i_1^1,u}^1\le (3^2-1) \underset{(I^{\kk}_1,K_1)\in\LL{2},(j_1,p_1)\in \LL{1}^1}\sum\Bigg[\mathbb{P}\Bigg(\underset{\left[c_{p_1}^1+\frac{j_1S}{u^{2/\alpha_1(c_{p_1}+1)}},c_{p_1}^1+\frac{(j_1+1)S}{u^{2/\alpha_1{c_{p_1}+1}}}\right]\times V_{I^{\kk}_1,K_1}}\sup X(u)>u, \\ &&\underset{\left[c_{p_1}^1+\frac{(j_1+1)S+\sqrt{S}}{u^{2/\alpha_1(c_{p_1}+1)}},c_{p_1}^1+\frac{(j_1+2)S}{u^{2/\alpha_1{c_{p_1}+1}}}\right]\times V_{I'^{\kk'}_1,K'_1}}\sup X(u)>u\Bigg)+\mathbb{P}\Bigg(\underset{\left[c_{p_1}^1+\frac{(j_1+1)S}{u^{2/\alpha_1(c_{p_1}+1)}},c_{p_1}^1+\frac{(j_1+1)S+\sqrt{S}}{u^{2/\alpha_1{c_{p_1}+1}}}\right]\times V_{I'^{\kk'}_1,K'_1}}\sup X(u)>u\Bigg)\Bigg]\\
&&\le \tilde{\mathbb{C}}\underset{(I^{\kk}_1,K_1)\in\LL{2},(j_1,p_1)\in \LL{1}^1}\sum \left(\mathbb{C}S^4e^{-\mathbb{C}_1 S^{\alpha_1/2}}+\prod_{i=1}^2\mH_{\alpha_i}[0,1](C_U)^{1/\alpha_i}S^{3/2}\right)\Psi(u)
\EQNY
for suitably chosen constant $\tilde{\mathbb{C}}$. Note that in the last formula $V_{I'^{\kk'}_1,K'_1}$ is one of the adjacent sets of $V_{I^{\kk}_1,K_1}$, and the number of it is at most $3^2-1$.
Using the same arguments we can obtain similar upper bounds for $\Sigma_{i_1^2,u}^2, \Sigma_{i_2^1,u}^2$ and $\Sigma_{i_2^2,u}^2$. Consequently, the same reasoning as \eqref{eqSig2} yields
\BQNY
\lim_{S\rw\IF}\lim_{u\rw\IF}\frac{\Sigma^2_u}{\tLu\Psi(u)}=0,
\EQNY
\cH{hence the claim follows.} \QED

{\bf Acknowledgement}: The authors kindly acknowledge partial
support from Swiss National Science Foundation Project 200021-1401633/1 and
%The research of V.I. Piterbarg is supported by Russian Foundation for Basic Research, Project 11-01-00050-a.
%kindly acknowledge partial support by
the project RARE -318984, a Marie Curie International Research Staff Exchange Scheme Fellowship within the 7th European Community Framework Programme.

\end{document}